\documentclass[12pt]{article}
\usepackage{amsmath, amsthm, amssymb, stmaryrd, fullpage}
\usepackage[all]{xy}

\newtheorem{theorem}{Theorem}[subsection]
\numberwithin{equation}{theorem}
\newtheorem{lemma}[theorem]{Lemma}

\newtheorem{prop}[theorem]{Proposition}
\theoremstyle{definition}
\newtheorem{defn}[theorem]{Definition}
\newtheorem{example}[theorem]{Example}
\newtheorem{remark}[theorem]{Remark}
\newtheorem{convention}[theorem]{Convention}

\def\FF{\mathbb{F}}
\def\RR{\mathbb{R}}
\def\ZZ{\mathbb{Z}}
\def\calO{\mathcal{O}}
\newcommand{\gothm}{\mathfrak{m}}
\newcommand{\gotho}{\mathfrak{o}}
\newcommand{\dlog}{d\log}

\def\alg{\mathrm{alg}}
\def\an{\mathrm{an}}
\def\con{\mathrm{con}}
\def\sep{\mathrm{sep}}
\def\perf{\mathrm{perf}}

\def\dual{\vee}
\def\fp{\pi^{-1}}
\def\del{\partial}

\def\Gperf{\Gamma^{\perf}}
\def\Galg{\Gamma^{\alg}}

\def\be{\mathbf{e}}
\def\bv{\mathbf{v}}
\def\bw{\mathbf{w}}

\def\Gancon{\Gamma_{\an,\con}}
\def\Gcon{\Gamma_{\con}}

\DeclareMathOperator{\Ext}{Ext}
\DeclareMathOperator{\Gal}{Gal}
\DeclareMathOperator{\Hom}{Hom}
\DeclareMathOperator{\naive}{naive}
\DeclareMathOperator{\rank}{rank}

\newcounter{fixmectr}

\begin{document}

\title{The $p$-adic local monodromy theorem for fake annuli}
\author{Kiran S. Kedlaya \\ Department of Mathematics \\ Massachusetts
Institute of Technology \\ 77 Massachusetts Avenue \\
Cambridge, MA 02139 \\
\texttt{kedlaya@math.mit.edu}}
\date{November 10, 2006}

\maketitle

\begin{abstract}
We establish a generalization of the $p$-adic local
monodromy theorem (of Andr\'e, Mebkhout, and the author) in which
differential equations on rigid analytic annuli are replaced by differential
equations on so-called ``fake annuli''. 
The latter correspond loosely to completions of
a Laurent polynomial ring with respect to a monomial valuation.
The result represents a step towards a higher-dimensional version of
the $p$-adic local monodromy theorem (the ``problem of semistable reduction'');
it can also be viewed as a novel presentation of the original
$p$-adic local monodromy theorem.
\end{abstract}

\tableofcontents

\section{Introduction}

This paper proves a generalization of the $p$-adic local monodromy
theorem of 
Andr\'e \cite{andre}, Mebkhout \cite{mebkhout}, and the present
author \cite{me-local}.
That theorem, originally conjectured by Crew \cite{crewfin} 
as an analogue in
rigid ($p$-adic) cohomology of Grothendieck's local monodromy theorem in
\'etale ($\ell$-adic) cohomology, asserts the quasi-unipotence of differential
modules with Frobenius structure on certain one-dimensional rigid
analytic annuli.

The $p$-adic local monodromy theorem has far-reaching consequences in
the theory of rigid cohomology, particularly for curves \cite{crewfin}.
However, although one can extend it to a relative form
\cite[Theorem~5.1.3]{me-finite} to obtain some higher-dimensional results,
for some applications one needs a version of the monodromy theorem which
is truly higher-dimensional.

This theorem takes an initial step towards producing such a higher-dimensional
monodromy theorem, by proving a generalization
in which the role of the annulus is replaced by a somewhat
mysterious space, called a \emph{fake annulus}, described by certain rings
of multivariate power series. When there is only one variable, the space
is a true annulus, so the result truly generalizes the original monodromy
theorem. Indeed, this paper has the side effect of giving an exposition
of the original theorem, albeit one somewhat encumbered by extra notation
needed for the fake case.

In the remainder of this introduction, we
explain further the context in which the $p$-adic local monodromy theorem 
arises, introduce and justify the fake analogue, and outline the structure
of the paper.

\subsection{Monodromy of $p$-adic differential equations}

Let $K$ be a field of characteristic zero
complete with respect to a nonarchimedean absolute
value, whose residue field $k$ has characteristic $p>0$.
Suppose we are given a rigid analytic annulus over $K$ and a differential
equation on the annulus, i.e., a module equipped with an integrable 
connection. We now wish to define the ``monodromy around the puncture''
of this connection, despite not having recourse to the analytic
continuation we would use in the analogous classical setting.
In particular, we would like to construct a representation of an
appropriate \'etale fundamental group, whose triviality or unipotence
amounts to the existence of a full set of horizontal sections or log-sections;
the latter relates closely
to the existence of an extension or logarithmic extension of the connection
across the puncture (e.g., \cite[Theorem~6.4.5]{me-part1}).

We can define a monodromy representation associated to a connection if 
we can find enough horizontal sections on some suitable covering
space. In particular,
we are mainly interested in connections which 
become unipotent, i.e., can be filtered
by submodules whose successive quotients are trivial for the connection,
on some cover of the annulus
which is ``finite \'etale near the boundary'' (in a sense that can be made
precise).
One can then construct a monodromy representation, using the Galois action
on horizontal sections, giving
an equivalence of categories between such \emph{quasi-unipotent}
modules with connection and a certain representation category.
(Beware that if the field $k$ is not algebraically closed, these
representations are only semilinear over the relevant field,
namely the maximal unramified extension of $K$. See 
\cite[Theorem~4.45]{me-mono-over} for a precise statement.)

In order for such an equivalence to be useful, we need to be able
to establish conditions under which a module with connection is
forced to be quasi-unipotent. As suggested in the introduction
to \cite{me-part1}, a natural, geometrically meaningful
candidate restriction (analogous to
the existence of a variation of Hodge structure for a complex
analytic connection) is the existence of a \emph{Frobenius structure}
on the connection. For $K$ discretely valued, the fact that
connections with a Frobenius structure are quasi-unipotent is the
content of the $p$-adic local monodromy theorem ($p$LMT) 
of Andr\'e \cite{andre}, Mebkhout \cite{mebkhout},
and this author \cite{me-local}.

Note that in this paper, we will not go all the way to the construction
of monodromy representations. These appear directly in Andr\'e's proof
of the $p$LMT (at least for $k$ algebraically closed); 
for a direct construction (of Fontaine type)
assuming the $p$LMT, see \cite{me-mono-over}.
See Remark~\ref{R:matsuda} for more details.

\subsection{Fake annuli}

The \emph{semistable reduction problem} (or \emph{global quasi-unipotence problem})
for overconvergent $F$-isocrystals,
as formulated by Shiho
\cite[Conjecture~3.1.8]{shiho2} and reformulated in
\cite[Conjecture~7.1.2]{me-part1}, is essentially to give a higher
dimensional version of the $p$LMT. From the point of view of \cite{me-part1},
this can be interpreted as proving a uniform version of the $p$LMT across
all divisorial valuations on the function field of the original variety.
This interpretation immediately suggests that one needs to exploit the
quasi-compactness of the Riemann-Zariski 
space associated to the function field of an irreducible variety; this
observation is developed in more detail in \cite{me-part2}.

The upshot is that one must prove the $p$LMT uniformly for the
divisorial valuations in a neighborhood (in Riemann-Zariski space) of an
arbitrary valuation, not just a divisorial one. Ideally, one could proceed
by first verifying whether the $p$LMT itself makes sense and
continues to hold true when one passes from a divisorial valuation to
a more general one. This entails replacing the annuli
in the $p$LMT with some sort of ``fake annuli'' which cannot be
described as rigid analytic spaces in the usual sense. 
Nonetheless,
one can still sensibly define rings of analytic functions in a neighborhood
of an irrational point, and thus set up a ring-theoretic framework in which
an analogue of the $p$-adic local monodromy theorem can be formulated.
(This allows us to get away with the linguistic swindle of speaking
meaningfully about ``$p$-adic differential equations on fake annuli'' without
giving the noun phrase ``fake annulus'' an independent meaning!)

Indeed, this framework fits naturally into the context of the slope
filtration theorem of \cite{me-local}. That theorem, which gives a
structural decomposition of a semilinear endomorphism on a finite free
module over the Robba ring (of germs of rigid analytic functions on an open
annulus with outer radius 1), does not make any
essential use of the fact that the Robba ring is described in terms of power
series. Indeed, as presented in \cite{me-slope}, the theorem applies directly
to our fake annuli; thus to prove the analogue of the $p$-adic monodromy 
theorem, one needs only analogize Tsuzuki's unit-root monodromy theorem
from \cite{tsuzuki-unitroot}. With a bit of effort, this can indeed be done,
thus illustrating some of the power of the Frobenius-based approach
to the monodromy theorem. Note that our definition of fake annuli
will actually include true rigid analytic annuli, so the monodromy
theorem given here will strictly generalize the $p$-adic local monodromy
theorem.

Unfortunately, it is not so clear how to prove a form of the $p$LMT for
arbitrary valuations; in this paper, we restrict to a somewhat simpler class.
These are the monomial valuations, which correspond
to monomial orderings in a polynomial ring. For instance, these include
valuations
on $k(x,y)$ in which the valuations of $x$ and $y$ are linearly independent
over the rational numbers. There are many
valuations that do not take this form, namely the \emph{infinitely singular}
valuations; the semistable reduction problem for
these must be treated in a more roundabout fashion, which we will not discuss
further here.

\subsection{Structure of the paper}

We conclude this introduction with a summary of the various
sections of the paper.

In Section~\ref{sec:fake}, we define the rings corresponding to fake
annuli, and verify that they fit into the formalism within which
slope filtrations are constructed in \cite{me-slope}.

In Section~\ref{sec:frob conn}, we define $F$-modules,
$\nabla$-modules, and $(F, \nabla)$-modules on fake annuli,
and verify that the category of $(F, \nabla)$-modules is invariant
of the choice of a Frobenius lift (Proposition~\ref{P:change Frob}).

In Section~\ref{sec:unit-root}, we give a fake annulus generalization
of Tsuzuki's theorem on unit-root $(F, \nabla)$-modules
(Theorem~\ref{T:tsuzuki}).

In Section~\ref{sec:monodromy}, we invoke the technology of slope filtrations
from \cite{me-local} (via \cite{me-slope}), and apply it to deduce from
Theorem~\ref{T:tsuzuki} a form of the $p$-adic
local monodromy theorem for $(F, \nabla)$-modules on
fake annuli (Theorem~\ref{T:local mono}).

In Section~\ref{sec:complements},
we deduce some consequences of the $p$-adic local monodromy theorem.
Namely, we calculate some extension groups in the category of
$(F, \nabla)$-modules, establish a local duality theorem,
and generalize some results from
\cite{dejong} and \cite{me-full} 
on the full faithfulness of overconvergent-to-convergent
restriction.

\subsection*{Acknowledgments}

Thanks to Nobuo Tsuzuki for some helpful remarks on the unit-root
local monodromy theorem.
Thanks also to Francesco Baldassarri and Pierre Berthelot for organizing
useful workshops on $F$-isocrystals and rigid cohomology in 
December 2004 and June 2005.
The author was supported by NSF grant number DMS-0400727.

\section{Fake annuli}
\label{sec:fake}

In this section, we describe ring-theoretically the fake annuli to
which we will be generalizing the $p$-adic local monodromy theorem,
deferring to \cite{me-slope} for most of the heavy lifting.

First, we put in some notational conventions that will hold in force
throughout the paper; for the most part, these hew to the notational
r\'egime of \cite{me-slope} (which in turn mostly follows
\cite{me-local}), with a few modifications made for greater consistency with
\cite{me-part1}.

\setcounter{equation}{0}
\begin{convention}
Throughout this paper, let $K$ be a complete
\emph{discretely valued} field of characteristic $0$, whose residue
field $k$ has characteristic $p>0$.
Let $\gotho = \gotho_K$ be the ring of integers
of $K$, and let $\pi$ denote a uniformizer of $K$.
Let $w$ be the valuation on $\gotho$ normalized so that
$w(\pi) = 1$.
Let $q$ be a power of $p$, and assume extant and fixed
a ring endomorphism $\sigma_K: K \to K$, continuous with respect to
the $\pi$-adic valuation, and lifting the $q$-power endomorphism on $k$.
Let $K_q$ be the fixed field of $K$ under $\sigma_K$,
and let $\gotho_q$ be the fixed ring of $\gotho$ under $\sigma_K$.
Finally,
let $I_n$ denote the $n \times n$ identity matrix over any ring.
\end{convention}

\begin{remark}
As noted in the introduction, the restriction to $K$ discretely
valued is endemic to the methods of this paper; see Remark~\ref{R:no
spherical} for further discussion.
\end{remark}

\subsection{Monomial fields}

\begin{defn} \label{D:monomial}
Let $k$ be a field.
A \emph{nearly monomial field (of height $1$) over $k$} 
is a field $E$ equipped with a valuation $v: E^* \to \RR$ (also written
$v: E \to \RR \cup \{+\infty\}$)
satisfying the following
restrictions.
\begin{enumerate}
\item[(a)]
The field $E$ is a separable extension of $k$, i.e., $k \subseteq E$
and $k \cap E^p = k^p$.
\item[(b)]
The image $v(E^*)$ of $v$ is a finitely generated $\ZZ$-submodule of $\RR$,
and $v(k^*) = \{0\}$.
\item[(c)]
The field $E$ is complete with respect to $v$.
\item[(d)]
With the notations
\begin{align*}
\gotho_E &= \{x \in E: v(x) \geq 0\} \\
\gothm_E &= \{x \in E: v(x) > 0\} \\
\kappa_E &= \gotho_E/\gothm_E,
\end{align*}
the natural map $k \to \kappa_E$ is finite.
\end{enumerate}
If $\kappa_E = k$, we say $E$ is a \emph{monomial field} 
(or \emph{fake power series field}) over $k$;
in that case, $k$ is integrally closed in $E$.
We define the \emph{rational rank} of $E$ to be the rank of $v(E^*)$ as
a $\ZZ$-module.
\end{defn}

\begin{remark}
One can also speak of monomial fields of height greater than 1, by allowing
the valuation $v$ to take values in a more general totally ordered 
abelian group. The techniques used in this paper do not apply to that case,
so we will ignore it; in the semistable reduction context, one can
eliminate the case of height greater than 1 by an inductive argument
\cite[Proposition~4.2.4]{me-part2}.
\end{remark}

\begin{example}
A monomial field of rational rank 1
is just a power series field, by the Cohen structure theorem.
This characterization generalizes to arbitrary monomial fields;
see Definition~\ref{D:coordinate}.
Note also that nearly monomial fields are examples of \emph{Abhyankar valuations},
i.e., valuations in which equality holds in Abhyankar's inequality
\cite[Th\'eor\`eme~9.2]{vaquie}.
\end{example}

\begin{remark}
If $E$ is a nearly monomial field and $E'/E$ is a finite separable extension, then
$E'$ is also nearly monomial: $E'$ is separable over $k$,
the valuation $v$ extends uniquely to a valuation $v'$ on $E'$,
$E'$ is complete with respect to $v'$, and the index $[v'((E')^*):v(E^*)]$ and
degree $[\kappa_{E'}:\kappa_E]$ are both finite since their product is at most
$\deg(E'/E)$ \cite[Proposition~5.1]{vaquie}. Conversely, every nearly monomial field
can be written as a finite separable extension of a monomial field; see
Definition~\ref{D:coordinate}.
\end{remark}

\begin{remark} \label{R:not fake again}
If $E$ is a nearly monomial field over $k$ and $\kappa_E/k$ is separable, then by Hensel's
lemma, the integral closure $k'$ of $k$ in $E$ is isomorphic to $\kappa_E$; in other words,
$E$ is a monomial field over $k'$. In particular, if $k$ is perfect, then any finite
extension of a monomial field over $k$ is a monomial field over some finite
separable extension of $k$. This fails if $k$ is not perfect, even for finite
separable extensions of the monomial field: the field $k$ is integrally
closed in
\[
k((t))[z]/(z^p - z - ct^{-p}) \qquad (c \in k \setminus k^p),
\]
but the latter has residue field $k(c^{1/p}) \neq k$.
\end{remark}

It will frequently be convenient to work with monomial fields
in terms of coordinate systems.
\begin{defn}
Let $m$ be a nonnegative integer.
Let $L$ be a lattice in $\RR^m$,
i.e., a $\ZZ$-submodule of $\RR^m$ which is free of rank $m$,
and which spans $\RR^m$ over $\RR$.
Let $L^\dual \subseteq (\RR^m)^\dual$ denote the lattice dual to $L$:
\[
L^\dual = \{\mu \in (\RR^m)^\dual: \mu(z) \in \ZZ \quad \forall z \in L\}.
\]
Given a formal sum $\sum_{z \in L} c_z \{z\}$, with the $c_z$ in some ring,
define the \emph{support} of the sum to be
the set of
$z \in L$ such that $c_z \neq 0$; define the support of a matrix of formal
sums to be the union of the supports of the entries.
If $S \subseteq L$ and a formal sum or matrix has support contained in $S$,
we also say that the element or matrix is ``supported on $S$''.
For $R$ a ring, 
let $R[L]$ denote the group algebra of $L$ over $R$, i.e., the set of
formal sums $\sum_{z \in L} c_z \{z\}$ with coefficients in $R$ and 
finite support.
\end{defn}

\begin{remark}
It is more typical to denote the class in $R[L]$ of a lattice element
$z \in L$ by $[z]$, rather than $\{z\}$.
However, we need to use brackets to denote Teichm\"uller lifts, so we
will stick to braces for internal consistency.
\end{remark}

\begin{defn}
For any ring $R$ and any $\lambda \in (\RR^n)^\dual$,
let $v_\lambda$ denote the valuation on $R[L]$ given by
\[
v_\lambda\left(\sum_{z \in L} c_z \{z\}\right) 
= \min\{\lambda(z): z \in L,\, c_z \neq 0\}.
\]
Let $P_\lambda \subseteq L$ denote the submonoid of $z \in L$ for
which $\lambda(z) 
\geq 0$, and let $R[L]_\lambda$ denote the monoid algebra $R[P_\lambda]$.
Let $R \llbracket L \rrbracket_\lambda$ and 
$R ((L))_\lambda$ denote the $v_\lambda$-adic completions
of $R[L]_\lambda$ and $R[L]$, respectively.
\end{defn}

\begin{defn} \label{D:coordinate}
Given a lattice $L$ and some $\lambda \in (\RR^m)^\dual$, we say
$\lambda$ is \emph{irrational} if 
$L \cap \ker(\lambda) = \{0\}$. In this case,
$k((L))_\lambda$ is a monomial field over $k$.
Conversely,
given a nearly monomial field $E$ over a field $k$, with valuation $v$, a
\emph{coordinate system} for $E$ is a sequence
$x_1, \dots, x_m$ of elements of $E$ such that
$v(x_1), \dots, v(x_m)$ freely generate $v(E^*)$ as a $\ZZ$-module.
Given a coordinate system, put $L = \ZZ^m$ with generators $z_1, \dots, z_m$,
and define $\lambda \in (\RR^m)^\dual$ by $\lambda(z_i) = v(x_i)$;
then $\lambda$ is irrational, and the
continuous map $k((L))_\lambda \to E$ given by $\{z_i\}
\mapsto x_i$ is injective. If we identify $k((L))_\lambda$ with its
image in $E$, then $E$ is finite separable over $k((L))_\lambda$; if $E$
is monomial over $k$, then in fact $E = k((L))_\lambda$ by 
Proposition~\ref{P:defectless} below.
This fact 
may be viewed as a
monomial version of the Cohen structure theorem in equal characteristics.
\end{defn}

\begin{prop} \label{P:defectless}
Let $L$ be a lattice in $\RR^m$, choose $\lambda \in (\RR^m)^\dual$
irrational, and let $E$ be a finite separable extension of $k((L))_\lambda$
with value group $\lambda(L)$ and residue field $k$. Then $E = k((L))_\lambda$.
\end{prop}
\begin{proof}
The claim is equivalent to showing that $k((L))_\lambda$ 
is inseparably defectless
in the sense of Ostrowski's lemma  \cite[Th\'eor\`eme~2, p.\ 236]{ribenboim}.
In particular, since a tamely ramified extension of $E$ is necessarily
without defect, it suffices to check that there is no Artin-Schreier defect
extension. If 
$E = 
k((L))_\lambda[z]/(z^p - z - x)$ were such an extension, we could rewrite it
as $k((L))_\lambda[z]/(z^p - z - y)$ with the leading term of $y$
being either an element of $pL$ times a non-$p$-th power $c$ in $k$, or an
element of $L \setminus pL$ times a nonzero element of $k$.
But in the first case the residue field of $E$ would be $k(c^{1/p})$,
and in the second case we would have
$[v(E):\lambda(L)] = p$; in either case, we would 
contradict the assumption that
$E$ is a defect extension. This contradiction yields the claim.
\end{proof}

\begin{remark}
The term ``monomial field'' is modeled on the use of the term ``monomial valuation'',
e.g., in \cite{favre-jonsson}, to refer to a valuation $v$ of the sort considered
in Definition~\ref{D:monomial}. (Such valuations, each of which endows the lattice
$L$ with a total ordering, are more common in mathematics than one
might initially realize: for example,
they are used to define highest weights in the theory of Lie algebras,
and they are sometimes used to construct term orders in the theory of Gr\"obner bases.)
In a previous version of this paper, the term ``fake power series field''
was used instead; we have decided that it would be better to save this term
for describing the completion of a finitely
generated field extension of $k$ with respect to \emph{any} valuation 
of height 1.
(See Remark~\ref{R:not rank 1} for some reasons why we are not considering 
such 
valuations here.)
\end{remark}

\subsection{Witt rings and Cohen rings}
\label{subsec:setup}

We now enter the formalism of \cite[\S~2]{me-slope}.
\begin{defn}
Let $K^{\perf}$ be the completion of the direct limit
$K \stackrel{\sigma_K}{\to} K \stackrel{\sigma_K}{\to} \cdots$ for the
$\pi$-adic topology; this is a complete discretely valued field of
characteristic 0 with residue field $k^{\perf}$, so it contains
$W(k^{\perf})$ by Witt vector functoriality.
For $E$ a perfect field of characteristic $p$ containing $k$, put
$\Gamma^E = W(E) \otimes_{W(k^{\perf})} \calO_{K^{\perf}}$; note that the
valuation $w$ extends naturally to $\Gamma^E$.
\end{defn}

\begin{defn} \label{D:partial valuation}
For $E$ a perfect field of characteristic $p$ containing $k$,
complete for a valuation $v$ trivial on $k$,
define the \emph{partial valuations} $v_n$ on $\Gamma^E[\fp]$ as follows.
Given $x \in \Gamma^E[\fp]$, write
$x = \sum_i [\overline{x_i}] \pi^i$, where each
$\overline{x_i} \in E$ and the brackets denote
Teichm\"uller lifts. Set
\[
v_n(x) = \min_{i \leq n} \{v(\overline{x_i})\}.
\]
As in \cite[Definition~2.1.5]{me-slope}, the partial
valuations satisfy some useful identities (here using
$\sigma$ to denote the $q$-power Frobenius):
\begin{align*}
v_n(x+y) &\geq \min\{v_n(x), v_n(y)\}  & (x,y \in \Gamma^E[\fp],\, 
n \in \ZZ) \\
v_n(xy) &\geq \min_{m \in \ZZ} \{v_m(x) + v_{n-m}(y)\}  & (x,y \in 
\Gamma^E[\fp],\, n \in \ZZ) \\
v_n(x^{\sigma}) &= q v_n(x) &(x \in \Gamma^E[\fp],\, n \in \ZZ) \\
v_n([\overline{x}]) &= v(\overline{x})  &(\overline{x}
 \in E,\, n \geq 0).
\end{align*}
In the first two cases, equality holds whenever the minimum is
achieved exactly once.
Define the \emph{levelwise topology} (or \emph{weak topology})
on $\Gamma^E$ by declaring that
a sequence $\{x_i\}$ converges to zero if and only if for each $n$,
$v_n(x_i) \to \infty$ as $i \to \infty$.
\end{defn}

\begin{defn}
For $r > 0$, write $v_{n,r}(x) = r v_n(x) + n$; for $r=0$,
write conventionally
\[
v_{n,0}(x) = \begin{cases} n & v_n(x) < \infty \\
\infty & v_n(x) = \infty.
\end{cases}
\]
Let $\Gamma^E_r$ be the subring of $\Gamma^E$
for which $v_{n,r}(x) \to \infty$ as $n \to \infty$; then
$\sigma$ sends $\Gamma^E_{r}$ to $\Gamma^E_{r/q}$.
Define the map $w_r$ on $\Gamma^E_r$ by
\[
w_r(x) = \min_n \{ v_{n,r}(x) \};
\]
then $w_r$ is a valuation on $\Gamma_r$ by
\cite[Lemma~2.1.7]{me-slope},
and $w_r(x) = w_{r/q}(x^{\sigma})$.
Put
\[
\Gamma^E_{\con} = \cup_{r>0} \Gamma^E_r;
\]
this is a henselian discrete valuation ring with maximal ideal
$\pi \Gamma^E_{\con}$ and residue field $E$
(see discussion in \cite[Definition~2.2.13]{me-slope}).
\end{defn}

\begin{convention}
For $E$ a not necessarily perfect field complete for a valuation
$v$ trivial on $k$, we write $E^{\perf}$ and $E^{\alg}$
for the \emph{completed} (with respect to $v$)
perfect and algebraic closures of $E$.
When $E$ is to be understood,
we abbreviate $\Gamma^{E^{\perf}}$ and $\Gamma^{E^{\alg}}$ to
$\Gperf$ and $\Galg$, respectively. Note that this is consistent
with the conventions of \cite{me-slope} but \emph{not} with those
of \cite{me-local}, where the use of these superscripts is taken not
to imply completion.
\end{convention}

Since we are interested in constructing $\Gamma^E$ for
$E$ a monomial field,
which is not perfect, we must do a bit more work, as in
\cite[\S~2.3]{me-slope}.
\begin{defn} \label{D:good pairs}
Let $E$ be a nearly monomial field over $k$ with valuation $v$.
Let $\Gamma^E$ be a complete discrete valuation ring of characteristic
$0$ containing $\gotho$ and having residue field $E$, such that
$\pi$ generates the maximal ideal of $\Gamma^E$.
Suppose that $\Gamma^E$ is equipped with a \emph{Frobenius lift}, i.e.,
a ring endomorphism $\sigma$ extending $\sigma_K$ on $\gotho_K$
and lifting the $q$-power Frobenius map on $E$.
We may then embed $\Gamma^E$ into $\Gperf$ by mapping $\Gamma^E$
into the first term of the direct system 
$\Gamma^E \stackrel{\sigma}{\to} \Gamma^E \stackrel{\sigma}{\to} \cdots$,
completing the direct system, 
and mapping the result into $\Gperf$ via Witt vector functoriality.
In particular, we may use this embedding to induce partial valuations
and a levelwise topology on $\Gamma^E$, taking care to remember that these
depend on the choice of $\sigma$. 
If $E'$ is a finite separable extension of $E$,
and we start with a suitable $\Gamma^E$ equipped with a Frobenius lift $\sigma$, we may
form the unramified extension of $\Gamma^E$ with residue field $E'$; this
will be a suitable $\Gamma^{E'}$, and carries a unique Frobenius lift extending
$\sigma$.
\end{defn}

\begin{defn}
Let $E$ be a nearly monomial field over $k$, and fix a pair
$(\Gamma^E, \sigma)$ as in Definition~\ref{D:good pairs}.
Write
\[
\Gamma^E_{\con} = \Gamma^E \cap \Gamma^{\perf}_{\con},
\]
with the intersection taking place within $\Gamma^{\perf}$.
For $r>0$, we say that $\Gamma^E$ \emph{has enough $r$-units}
if $\Gamma^E \cap \Gamma^{\perf}_r$ 
contains units lifting all nonzero elements
of $E$. We say that $\Gamma^E$ \emph{has enough units} (or more properly,
the pair $(\Gamma^E, \sigma)$ has enough units) if $\Gamma^E$ has
enough $r$-units for some $r>0$; this implies that $\Gamma^E_{\con}$
is a henselian discrete valuation ring with maximal ideal 
$\pi \Gamma^E_{\con}$ and residue field $E$.
If $\Gamma^E$ has enough units, then so does
$\Gamma^{E'}$ for any finite separable extension $E'$ of $E$
\cite[Lemma~2.2.12]{me-slope}. 
\end{defn}

\subsection{Toroidal interpretation}

The condition of having enough units is useful in the theory of slope
filtrations, but is not convenient to check in practice. Fortunately,
it has a more explicit interpretation in terms of certain
``na\"\i ve'' analogues of the functions $v_n$ and $w_r$,
as in \cite[\S~2]{me-local} or \cite[\S~2.3]{me-slope}.

\begin{defn}
Let $L$ be a lattice in $\RR^m$ and let $\lambda \in (\RR^m)^\dual$
be an irrational linear functional. Let $\Gamma^\lambda$ denote the
$\pi$-adic completion of $\gotho ((L))_\lambda$; its elements may be
viewed as formal sums $\sum_{z \in L} c_z\{z\}$ with $w(c_z) \to \infty$
as $\lambda(z) \to -\infty$. Define the \emph{na\"\i ve partial valuations}
on $\Gamma^\lambda[\fp]$ by the formula
\[
v^{\naive}_n\left(\sum c_z \{z\}\right) 
= \min\{ \lambda(z): z\in L,\, w(c_z) \leq n\},
\]
where the minimum is infinite if the set of candidate $z$'s is empty.
These functions satisfy the identities
\begin{align*}
v^{\naive}_n(x+y) &\geq \min\{v^{\naive}_n(x), v^{\naive}_n(y)\} 
& (x,y \in \Gamma^\lambda[\fp])\\
v^{\naive}_n(xy) &\geq \min_{m \in \ZZ} \{ v^{\naive}_m(x) + 
v^{\naive}_{m-n}(y)\} & (x,y \in \Gamma^\lambda[\fp]),
\end{align*}
with equality in each case if the minimum is achieved only once.
Define the \emph{na\"\i ve levelwise topology} 
(or \emph{na\"\i ve weak topology}) on $\Gamma^\lambda$ by declaring that
a sequence $\{x_i\}$ converges to zero if and only if for each $n$,
$v^{\naive}_n(x_i) \to \infty$ as $i \to \infty$.
\end{defn}

\begin{defn}
For $r > 0$ and $n \in \ZZ$, write
\[
v^{\naive}_{n,r}(x) = r v^{\naive}_n(x) + n;
\]
extend the definition to $r=0$ by setting
\[
v^{\naive}_{n,0}(x) = \begin{cases} n & v^{\naive}_n(x) < \infty \\
\infty & v^{\naive}_n(x) = \infty.
\end{cases}
\]
Let $\Gamma_r^{\naive}$ be the set of $x \in \Gamma^\lambda$ such that
$v^{\naive}_{n,r}(x) \to \infty$ as $n \to \infty$.
Define the map $w^{\naive}_r$ on $\Gamma_r^{\naive}$ by
\[
w^{\naive}_r(x) = \min_{n \in \ZZ} \{v^{\naive}_{n,r}(x)\};
\]
then $w^{\naive}_r$ is a valuation on $\Gamma_r^{\naive}[\fp]$,
as in \cite[Lemma~2.1.7]{me-slope}.
Put 
\[
\Gcon^{\naive} = \cup_{r>0} \Gamma_r^{\naive}.
\]
\end{defn}

\begin{remark}
The ring $\Gamma_r^{\naive}$ is a principal ideal domain;
this will follow from \cite[Proposition~2.6.5]{me-slope} in conjunction
with Definition~\ref{D:standard ext} below.
\end{remark}

We may view $\Gamma^\lambda$ as an instance of the definition of
$\Gamma^E$ in the case
$E = k((L))_\lambda$; this gives sense to the following result.
\begin{prop} \label{P:naive compare}
Let $\sigma$ be a Frobenius lift on $\Gamma^E = \Gamma^\lambda$
for $E = k((L))_\lambda$. Then for $r>0$, the following are
equivalent.
\begin{enumerate}
\item[(a)] $\sigma$ is continuous for the na\"\i ve levelwise topology
(i.e., that topology coincides with the levelwise topology induced
by $\sigma$),
and for each $z \in L$ nonzero,
$\{z\}^\sigma/\{z\}^q$ is a unit in $\Gamma_r^{\naive}$.
\item[(b)] For $s \in (0,qr]$, $n \geq 0$, and $x \in \Gamma^E$,
\begin{equation} \label{eq:compare naive}
\min_{j \leq n} \{v_{j,s}(x)\} = \min_{j \leq n} \{ v_{j,s}^{\naive}(x)
\}.
\end{equation}
\item[(c)] $\Gamma^E$ has enough $qr$-units,
and for each $z \in L$ nonzero, $\{z\}$ is a unit in $\Gamma^E_{qr}$.
\end{enumerate}
In particular, in each of these cases,
for $s \in (0,qr]$, $\Gamma^{\naive}_s = \Gamma^E_s$ 
and $w_s(x) = w_s^{\naive}(x)$ for all $x \in \Gamma_s$.
\end{prop}
\begin{proof}
Given (a), for $s \in (0,qr]$, we have
\begin{equation} \label{eq:compare naive1}
\min_{j \leq n} \{v_{j,s}^{\naive}(x) \} = 
\min_{j \leq n} \{v_{j,s/q}^{\naive}(x^\sigma) \}
\end{equation}
for each $n \geq 0$ 
and each $x \in \Gamma$, as in the proof of
\cite[Lemma~2.3.3]{me-slope}.
We then obtain (b) as in the proof of
\cite[Lemma~2.3.5]{me-slope}, from which (c) follows immediately.

Given (c), the equation \eqref{eq:compare naive} holds for 
$x = \{z\}$ for any $z \in L$, since the minima both occur for
$j=0$. For $x = \sum c_z \{z\}$ a finite sum, we have
\begin{equation} \label{eq:compare naive2}
\min_{j \leq n} \{v_{j,s}^{\naive}(x)\} = \min_{j \leq n}
\{ \min_{z \in L} \{v_{j,s}^{\naive}(c_z\{z\})\} \}
\end{equation}
and so the left side of \eqref{eq:compare naive} is greater than
or equal to the right side. 
On the other hand, if $j$ is taken to be the smallest
value for which the outer minimum is achieved on the right side of
\eqref{eq:compare naive2}, then the inner minimum is achieved by a
unique value of $z$. Thus we actually may deduce equality in
\eqref{eq:compare naive}, again for $x = \sum c_z \{z\}$ a finite sum.
For general $x$, we may obtain the desired equality by replacing $x$ 
by a finite sum $x'$ such that
$x-x' = y+z$ for some $y \in \Gamma^E$ with $w(y)$ greater than $n$,
and some $z \in \Gamma^E_{qr}$ with $w_{qr}(z)$ greater than either
side of \eqref{eq:compare naive}. Hence (c) implies (b).

Finally, note that (b) implies (a) straightforwardly.
\end{proof}
\begin{remark}
Note that in Proposition~\ref{P:naive compare}, conditions
(a) and (c) may be checked for $z$ running over a basis of $L$.
Note also
that Proposition~\ref{P:naive compare} implies that for $E = k((L))_\lambda$,
if $\Gamma^E$ has enough units, then $\Gamma^E$
is isomorphic to $\Gamma^\lambda$.
\end{remark}

\begin{defn} \label{D:standard ext}
By the \emph{standard extension} of $\sigma_K$ to
$\Gamma^\lambda$, we will mean the
Frobenius lift $\sigma$ defined by
\[
\sum_{z \in L} c_z \{z\} \mapsto \sum_{z \in L} c_z^{\sigma_K} \{z\}^q.
\]
(We will also refer to such a $\sigma$ as a \emph{standard Frobenius lift}.)
When equipped with a standard Frobenius lift, $\Gamma^\lambda$ 
has enough $r$-units
for every $r>0$; by Proposition~\ref{P:naive compare}, it follows that
$v_n(x) = v_n^{\naive}(x)$ 
for all $n \in \ZZ$ and all $x \in \Gamma^\lambda[\fp]$.
Thus many of the results of \cite[\S~2]{me-slope}, proved in terms
of the Frobenius-based valuations, also apply verbatim to the na\"\i ve
valuations.
\end{defn}

\begin{remark} \label{R:not rank 1}
For applications to semistable reduction,
one would also like to consider a similar situation in which
the residue field $k((L))_\lambda$ is
replaced by the completion of a finitely generated field extension of $k$
with respect to an arbitrary valuation of height (real rank) 1,
at least in the case where the transcendence degree over
$k$ is equal to 2.
This would require a slightly more flexible set of foundations: one must
work only with finitely generated $k$-subalgebras of the complete field,
so that one has hope of having enough units.
A more serious problem is
how to perform Tsuzuki's method (a/k/a Theorem~\ref{T:tsuzuki}) in this
context.
\end{remark}

\subsection{Analytic rings}

We now introduce ``analytic rings'', citing into \cite{me-slope}
for their structural properties.

\begin{defn}
Let $E$ be a nearly monomial field over $k$, or the completed perfect or algebraic
closure thereof.
In the first case, suppose that $\Gamma^E$ has enough $r_0$-units for some $r_0 > 0$
(otherwise take $r_0 = \infty$).
Let $I$ be a subinterval of
$[0,r_0)$ bounded away from $r_0$ (i.e., $I$ is a subinterval of $[0,r]$
for some $r \in (0,r_0)$).
Let $\Gamma^E_I$ denote the Fr\'echet completion of $\Gamma^E_{r_0}[\fp]$
under the valuations $w_s$ for $s \in I$;
this ring is an integral domain \cite[Lemma~2.4.6]{me-slope}.
If $I$ is closed, then
$\Gamma^E_I$ is a principal ideal domain \cite[Proposition~2.6.9]{me-slope}.
Put
\[\Gamma^E_{\an,r} = \Gamma^E_{(0,r]};
\]
this ring is a B\'ezout ring,
i.e., a ring in which
 every finitely generated ideal is principal
\cite[Theorem~2.9.6]{me-slope}.
Put $\Gancon^E = \cup_{r>0} \Gamma^E_{\an,r}$; then $\Gancon^E$
is also a B\'ezout ring.
The group of units in $\Gancon^E$ consists of the
nonzero elements of $\Gcon[\fp]$ \cite[Corollary~2.5.12]{me-slope}.
For $E'$ finite separable over $E$, $E' = E^{\perf}$, or
$E' = E^{\alg}$, by \cite[Proposition~2.4.10]{me-slope} one has
\[
\Gancon^{E'} = \Gancon^E \otimes_{\Gcon^E} \Gcon^{E'}.
\]
\end{defn}

\begin{remark}
It is likely that $\Gamma^E_I$ is a B\'ezout ring for any $I$
as above. However, this statement
is not verified in \cite{me-slope}, and we will not need it anyway,
so we withhold further comment on it.
\end{remark}

\begin{remark} \label{R:true}
If $E = k((t))$ is a power series field, then
the ring $\Gamma^E_I$ is the ring of rigid analytic functions on the
annulus $w(t) \in I$ in the $t$-plane.
Thus our construction
of fake annuli includes ``true'' one-dimensional
rigid analytic annuli over $K$, and most
of our results on fake annuli
(like the $p$-adic local monodromy theorem) generalize extant
theorems on true annuli. 
On the other hand, if $E = k((L))_\lambda$ and $\rank(L) > 1$,
then the ring $\Gamma^E_I$ is trying
to be the ring of rigid analytic functions on a subspace of the rigid
affine plane in the variables $\{z_1\}, \dots, \{z_m\}$ for some
basis $z_1, \dots, z_m$ of $L$, consisting of points for which
there exists $r \in I$ with $w(\{z_i\}) = r \lambda(z_i)$
for $i=1, \dots, m$.
If $I = [r,r]$, then this space is an affinoid space in the sense of
Berkovich,
but otherwise it is not (because one can only cut out an
analytic subspace of the form $w(x) = \alpha w(y)$ for $\alpha$ rational).
Indeed, as far as we can tell, this space is not a $p$-adic analytic space
in either of the Tate or Berkovich senses,
despite the fact that it has a sensible ring
of analytic functions; hence the use of the adjective ``fake'' in the
phrase ``fake annulus'', and the absence of an honest definition of that
phrase.
\end{remark}

\begin{remark} \label{R:no spherical}
Since one can sensibly define rigid analytic annuli over arbitrary complete
nonarchimedean fields, Remark~\ref{R:true}
suggests the possibility of working
with fake annuli over more general complete $K$.
However, the algebraic issues here get more complicated,
and we have not straightened them out to our satisfaction. 
For example, the analogue of the ring $\Gamma_r^{\naive}$ fails to be
a principal ideal ring if the valuation on $K$ is not discrete; it probably
still has the B\'ezout property (that finitely generated ideals are
principal), but we have not checked this.
In any case, the formalism of \cite{me-slope} completely breaks down when $K$
is not discretely valued, so an attempt here to avoid a discreteness hypothesis
now would fail to improve upon our ultimate results; we have thus
refrained from making such an attempt.
\end{remark}

\section{Frobenius and connection structures}
\label{sec:frob conn}

We now introduce a notion which should be thought of as 
a $p$-adic differential equation
with Frobenius structure on a fake annulus.
We start with some notational conventions.

\setcounter{theorem}{0}
\begin{convention} \label{conv:fake power}
Throughout this section, assume that $E$ is a monomial field
and that $\Gamma^E$ is equipped with a Frobenius lift such that
$\Gamma^E$ has enough $r_0$-units for some $r_0 > 0$; we view 
$\Gamma^E$ as being equipped with a levelwise topology via the choice
of a coordinate system. (This choice does not matter, as the topology can
be characterized as the coarsest one under which the $v_{n,r}$ are
continuous for all $n \in \ZZ$ and all $r \in (0,r_0)$.) We suppress
$E$ from the notation, writing $\Gamma$ for $\Gamma^E$,
$\Gcon$ for $\Gcon^E$, and so on.
\end{convention}

\begin{convention}
When a valuation is applied to a matrix, it is defined to be the minimum
value over entries of the matrix.
\end{convention}

We also make a definition of convenience.
\begin{defn} \label{D:admissible}
Under Convention~\ref{conv:fake power}, we will mean by an
\emph{admissible ring} any one of the following topological rings.
\begin{itemize}
\item
The ring $\Gamma$ or $\Gamma[\fp]$ with its levelwise topology.
\item
The ring $\Gamma_r$ or $\Gamma_r[\fp]$ with the Fr\'echet 
topology induced by $w_s$
for all $s \in (0,r]$, for $r \in (0,r_0)$. Note that for $\Gamma_r$,
this coincides with the topology induced by $w_r$ alone.
\item
The ring
$\Gcon$ or $\Gcon[\fp]$ topologized as the direct limit of the $\Gamma_r$
or $\Gamma_r[\fp]$.
\item
The ring $\Gamma_I$ with the Fr\'echet topology induced by the
$w_s$ for $s \in I$, for some $I \subseteq [0,r_0)$ bounded away from
$r_0$.
\item
The ring $\Gancon$ topologized as the direct limit of the $\Gamma_{\an,r}$.
\end{itemize}
By a \emph{nearly admissible ring}, we mean one of the above rings
with $E$ replaced by a finite separable extension.
\end{defn}

\subsection{Differentials}

\begin{defn}
Let $S/R$ be an extension of topological rings. A \emph{module of
continuous differentials} is a topological $S$-module $\Omega^1_{S/R}$
equipped with a continuous $R$-linear derivation $d: S \to \Omega^1_{S/R}$,
having the following universal property: for any topological $S$-module
$M$ equipped with a continuous $R$-linear derivation $D: S \to M$,
there exists a unique morphism $\phi: \Omega^1_{S/R} \to M$ of topological
$S$-modules such that $D = \phi \circ d$. Since the definition is via
a universal property, the module of continuous differentials is unique
up to unique isomorphism if it exists at all.
\end{defn}
Constructing modules of continuous differentials is tricky in general
(imitating the usual construction of the module of K\"ahler differentials
requires a topological tensor product, which is a rather delicate matter);
however, for fake annuli, the construction is straightforward.

\begin{defn}
By a \emph{coordinate system} for $\Gamma$, we will mean a
lattice $L$ in some $\RR^m$, an irrational linear functional
$\lambda \in (\RR^m)^\dual$,
an isomorphism $\Gamma^\lambda \cong \Gamma$ carrying
$z \in L$ to a unit in $\Gamma_{r_0}$ for each
nonzero $z \in L$, and
a basis $z_1, \dots, z_m$ of $L$. Such data always exist thanks
to Proposition~\ref{P:naive compare}.
\end{defn}

\begin{defn}
For the remainder of this subsection, choose 
a coordinate system for $\Gamma$, and let
$\mu_1, \dots, \mu_m \in L^\dual$ denote the basis dual to
$z_1, \dots, z_m$.
For $\mu \in L^\dual$ and $S$ an admissible ring, 
let $\del_\mu$ be the continuous derivation on $S$
defined by the formula
\[
\del_\mu\left(\sum_z c_z \{z\}\right) = \sum_z \mu(z) c_z \{z\};
\]
note that $\mu(z) \in \ZZ$, so it may sensibly be viewed as an element
of $\gotho$. (The continuity of $\del_\mu$ is clear in terms of na\"\i ve
partial valuations, so Proposition~\ref{P:naive compare} implies continuity
in terms of the Frobenius-based valuations.)
For $\mu = \mu_i$, write $\del_i$ for $\del_{\mu_i}$.
Define $\Omega^1_{S/\gotho}$ to be the free $S$-module
$S\,d\{z_1\} \oplus \cdots \oplus S\,d\{z_m\}$, equipped with the natural
induced topology and with the continuous $\gotho$-linear 
derivation $d: S \to \Omega^1_{S/\gotho}$ given by
\[
dx = \sum_{i=1}^m \del_{i}(x)\,\dlog \{z_i\}
\]
(where $\dlog(f) = df/f$).
\end{defn}
\begin{prop} \label{P:module of diff}
The module $\Omega^1_{S/\gotho}$ is a module of continuous derivations
for $S$ over $\gotho$. In particular, the construction does not depend on
the choice of the coordinate system.
\end{prop}
\begin{proof}
This is a straightforward consequence of the fact that
one of $\gotho[\{z_i\}^{\pm 1}]$ or
$\gotho[\fp, \{z_i\}^{\pm 1}]$ is dense in $S$.
\end{proof}

\begin{remark}
Note that Proposition~\ref{P:module of diff} also allows us to construct
the module of continuous differentials
$\Omega^1_{S/\gotho}$ when $S$ is only nearly admissible.
\end{remark}

\begin{remark} \label{R:cannot integrate}
For $\rank(L) = 1$ and $\mu \in L$ nonzero, the image of $\del_{\mu}$
is closed; however, this fails for $\rank(L)>1$,
because bounding $\lambda(z)$ does not in any way limit the
$p$-adic divisibility of $z$ within $L$.
This creates a striking difference between the milieux of true and fake annuli,
from the point of view of the study of differential equations.
On true annuli, one has the rich theory of $p$-adic differential
equations due to Dwork-Robba, Christol-Mebkhout, et al.
On fake annuli, much of that theory falls apart; the parts that
survive are those that rest upon Frobenius structures, whose behavior
differs little in the two settings.
\end{remark}

\subsection{$\nabla$-modules}

\begin{defn} 
Let $S$ be a nearly admissible ring.
Define a \emph{$\nabla$-module}
over $S$ to be a finite free $S$-module $M$ equipped with
an integrable $\gotho$-linear connection $\nabla: M \to M \otimes
\Omega^1_{S/\gotho}$; here integrability means that, letting
$\nabla_1$ denote the induced map 
\[
M \otimes_S \Omega^1_{S/\gotho} \stackrel{\nabla \otimes 1}{\to}
M \otimes_S \Omega^1_{S/\gotho} \otimes_S \Omega^1_{S/\gotho}
\stackrel{1 \otimes \wedge}{\to} M \otimes_S \wedge^2_S \Omega^1_{S/\gotho},
\]
the composite map $\nabla_1 \circ \nabla$ is zero.
We say $\bv \in M$ is \emph{horizontal} if $\nabla (\bv) = 0$.
\end{defn}

\begin{defn}
Suppose $S$ is admissible, and fix a coordinate system for $\Gamma$.
Given a $\nabla$-module $M$ over $S$, for $\mu \in L^\dual$,
define the map $\Delta_\mu: M \to M$ by writing
$\nabla(\bv) = \sum_{i=1}^m \bw_i \otimes \dlog\{z_i\}$ with
$\bw_i \in M$, and setting
\[
\Delta_\mu(\bv) = \sum_{i=1}^m \mu(z_i) \bw_i.
\]
Also, write $\Delta_i$ for $\Delta_{\mu_i}$.
\end{defn}

\begin{remark}
The maps $\Delta_\mu$ satisfy the following properties.
\begin{itemize}
\item
The map $L^\dual \times M \to M$ given by $(\mu, \bv) \mapsto \Delta_\mu(\bv)$
is additive in each factor.
\item
For all $\mu \in L^\dual$, $s \in S$, and $\bv \in M$, we have
the Leibniz rule
\[
\Delta_\mu(s\bv) = s \Delta_\mu(\bv) + \del_\mu(s) \bv.
\]
\item
For $\mu_1, \mu_2 \in L^\dual$, the maps $\Delta_{\mu_1}, \Delta_{\mu_2}$ commute.
\end{itemize}
Conversely, given a finite free $S$-module $M$ equipped with
maps $\Delta_\mu: M \to M$ for each $\mu \in L^\dual$ satisfying these
conditions, one can uniquely reconstruct a $\nabla$-module structure on $M$
that gives rise to the $\Delta_\mu$.
\end{remark}

\begin{remark}
Note that for true annuli (i.e., $\rank(L)=1$), the integrability restriction
is empty because $\Omega^1_{S/\gotho}$ has rank 1 over $S$. However,
for fake annuli, integrability is a real restriction: even though
the ring theory looks one-dimensional, the underlying ``fake space'' is
really $m$-dimensional, inasmuch as $\Omega^1_{S/\gotho}$ has rank $m$
over $S$.
\end{remark}

\begin{defn}
Let $M$ be a $\nabla$-module over $\Gancon^{E'}$,
for $E'$ a finite separable extension of $E$. We say $M$ is:
\begin{itemize}
\item
\emph{constant} if $M$ admits
a horizontal basis (a basis of elements of the kernel of $\nabla$);
\item
\emph{quasi-constant} if there exists a finite separable extension $E''$
of $E'$ such that $M \otimes \Gancon^{E''}$ is constant;
\item
\emph{unipotent} if $M$ admits an exhaustive filtration by 
saturated $\nabla$-submodules, whose successive quotients are
constant;
\item
\emph{quasi-unipotent} if $M$ admits an exhaustive filtration by
saturated $\nabla$-submodules, whose successive quotients are 
quasi-constant.
\end{itemize}
We extend these definitions to $(F, \nabla)$-modules by applying them
to the underlying $\nabla$-module.
\end{defn}
\begin{remark}
If $M$ is quasi-unipotent,
then there exists a finite separable extension $E''$ of $E'$
such that $M \otimes \Gancon^{E''}$ is unipotent.
The converse is also true: if $M \otimes \Gancon^{E''}$ is unipotent,
then the shortest unipotent filtration of $M \otimes \Gancon^{E''}$
is unique, so descends to $\Gancon^{E'}$.
\end{remark}

\subsection{Frobenius structures}

\begin{defn}
Let $S$ be a nearly admissible ring stable under $\sigma$;
for instance, $\Gamma, \Gcon, \Gancon$ are permitted, but
$\Gamma_r$ is not.
Define an \emph{$F$-module} over $S$ (with respect to $\sigma$)
to be a finite free $S$-module $M$ equipped with a $S$-module
homomorphism $F: \sigma^* M \to M$ which is an isogeny, i.e., which
becomes an isomorphism upon tensoring with $S[\fp]$. We typically
view $F$ as a $\sigma$-linear map from $M$ to itself;
we occasionally view $M$ as a left module for the twisted polynomial
ring $S\{\sigma\}$.
Given an $F$-module $M$ over $S$ and an integer $c$, which must be
nonnegative if $\pi^{-1} \notin S$, define the \emph{twist} $M(c)$ of $M$
to be a copy of $M$ with the action of $F$ multiplied by $\pi^c$.
\end{defn}

\begin{defn}
Let $S$ be a nearly admissible ring stable under $\sigma$.
Define an \emph{$(F,\nabla)$-module} over $S$
to be a finite free $S$-module $M$ equipped with the structures of
both an $F$-module and a $\nabla$-module, which are compatible in the
sense of making the following diagram commute:
\[
\xymatrix{
M \ar^-{\nabla}[r] \ar^{F}[d] & M \otimes \Omega^1_{S/\gotho} \ar^{F \otimes 
d\sigma}[d] \\ M \ar^-{\nabla}[r] & M \otimes \Omega^1_{S/\gotho}.
}
\]
\end{defn}

\begin{remark} \label{R:omega-sigma}
We may regard $\Omega^1_{S/\gotho}$ itself as an $F$-module via
$d\sigma$, in which case the compatibility condition asserts that
$\nabla: M \to M \otimes \Omega^1_{S/\gotho}$ is an $F$-equivariant map.
The fact that $\del_{\mu}(f^\sigma) \equiv 0 \pmod{\pi}$ for any
$f \in \Gcon$ means that $\Omega^1_{\Gcon/\gotho}$ 
is isomorphic as an $F$-module
to $N(1)$, for some $F$-module $N$ over $\Gcon$.
In the language of \cite{me-slope}, this means that the 
generic HN slopes of $\Omega^1_{\Gcon/\gotho}$ are positive
 \cite[Proposition~5.1.3]{me-slope}.
\end{remark}

\begin{defn}
For $a$ a positive integer, define an \emph{$F^a$-module} or
\emph{$(F^a, \nabla)$-module} as an $F$-module or $(F, \nabla)$-module
relative to $\sigma^a$. Given an $F$-module $M$, viewed as a
left $S\{\sigma\}$-module, define the $F^a$-module $[a]_* M$ to be
the left $S\{\sigma^a\}$-module given by restriction along the inclusion
$S\{\sigma^a\} \hookrightarrow S\{\sigma\}$; in other words, replace
the Frobenius action by its $a$-th power. Given an $F^a$-module $N$,
viewed as a left $S\{\sigma^a\}$-module, define the $F$-module
$[a]^* M$ to be the left $S\{\sigma\}$-module
\[
[a]^* M = S\{\sigma\} \otimes_{S\{\sigma^a\}} M;
\]
then
the functors $[a]^*$ and $[a]_*$ are left and right adjoints of each other.
See \cite[\S~3.2]{me-slope} for more on these operations.
\end{defn}

\subsection{Change of Frobenius}

The category of $(F, \nabla)$-modules over $\Gancon$ relative to $\sigma$
turns out to be canonically
independent of the choice of $\sigma$, by a Taylor series argument
(as in
\cite[\S~3.4]{tsuzuki-slope}).

\begin{convention}
Throughout this subsection,
fix a coordinate system on $\Gamma$.
Given an $m$-tuple $J = (j_1, \dots, j_m)$ of nonnegative integers,
write $J! = j_1! \cdots j_m!$; if $U = (u_1, \dots, u_m)$, write
$U^J = u_1^{j_1} \cdots u_m^{j_m}$. Also, define the ``falling factorials''
\begin{gather*}
\del^{\underline{J}} = \prod_{i=1}^m \prod_{l=0}^{j_i-1} (\del_i - l) \\
\Delta^{\underline{J}} = \prod_{i=1}^m \prod_{l=0}^{j_i-1} (\Delta_i - l),
\end{gather*}
with the convention that $\del^{\underline{0}}$ and $\Delta^{\underline{0}}$
are the respective identity maps.
(The use of falling factorial notation is modeled on \cite{gkp}.)
\end{convention}

\begin{lemma} \label{L:Leibniz rule}
Let $M$ be a $\nabla$-module over $\Gancon$. Then for any
$r \in \Gancon$, any $\bv \in M$, and any $m$-tuple $J$ of 
nonnegative integers,
\[
\frac{1}{J!} \Delta^{\underline{J}}(r\bv)
= \sum_{J_1+J_2=J} 
\left(\frac{1}{J_1!} \del^{\underline{J_1}}(r)\right)
\left(\frac{1}{J_2!} \Delta^{\underline{J_2}}(\bv) \right).
\]
\end{lemma}
\begin{proof}
Since
\[
\del_i(\del_i-1)\cdots(\del_i-j+1) = 
\{z_i\}^j (\{z_i\}^{-1} \del_i)^j
\]
and similarly for $\Delta_i$,
this amounts to a straightforward application of the Leibniz rule.
\end{proof}

\begin{lemma} \label{L:convergence1}
For any $u_1, \dots, u_m \in \Gcon$ with $w(u_i) > 0$ for
$i=1,\dots, m$, and any $x \in \Gancon$,
the series
\[
f(x) = \sum_{j_1, \dots, j_m =0}^\infty 
\frac{1}{J!} U^J \del^{\underline{J}}(x)
\]
converges in $\Gancon$, and the map $x \mapsto f(x)$ is a continuous
ring homomorphism sending $\{z_i\}$ to $u_i$.
\end{lemma}
\begin{proof}
Pick $r>0$ such that $u_1, \dots, u_m, x \in \Gamma_{\an,r}$
and $w_r(u_i) > 0$ for $i=1, \dots, m$.
Write $x = \sum_{z \in L} c_z\{z\}$; note that
for each $J$,
\[
\frac{1}{J!} \del^{\underline{J}}(x) = \sum_{z \in L} \left( \prod_{i=1}^m 
\binom{\mu_i(z)}{j_i} \right) c_z \{z\},
\]
so that $w_s(\del^{\underline{J}}(x)/J!) \geq w_s(x)$ for $s \in (0,r]$.
This yields the desired convergence, as well as continuity of the
map $x \mapsto f(x)$.
Moreover, $f$ is a ring homomorphism on
$\gotho[\{z_1\}, \dots, \{z_m\}]$ by Lemma~\ref{L:Leibniz rule},
so must be a ring homomorphism on $\Gancon$ by continuity;
the fact that it sends $\{z_i\}$ to $u_i$ is apparent from the formula.
\end{proof}

\begin{lemma} \label{L:bound nabla}
Let $M$ be a $\nabla$-module over $\Gamma_{\an,r}$ for some $r>0$.
Suppose that for some positive integer $h$,
$M$ admits a basis $\be_1, \dots, \be_n$ such that
the $n \times n$ matrices $N_1, \dots, N_m$ defined by
$\Delta_i (\be_l) = \sum_j (N_i)_{jl} \be_j$ satisfy
$w_r(N_i) > w((p^h)!)$ for $i=1, \dots, m$. 
For $J$ an $m$-tuple of nonnegative
integers, define the $n \times n$ matrix $N_J$ by
\[
\Delta^{\underline{J}} (\be_l) = \sum_j (N_J)_{jl} \be_j.
\]
Then
\[
w_r(N_J) \geq w(J!)-w(p)(j_1+\cdots+j_m)/(p^h(p-1)).
\]
\end{lemma}
\begin{proof}
The condition that $w_r(N_i) > w((p^h)!)$ means that
for any $a \in \ZZ$ and any $b \in \{0, \dots, p^h-1\}$, if we write
\begin{align*}
\bv &= \sum_{j=1}^m x_j \be_j & (x_j \in \Gamma_{\an,r}) \\
(\Delta_i-ap^h)(\Delta_i-ap^h-1)\cdots
(\Delta_i-ap^h-b)\bv &= \sum_{j=1}^m y_{ijab} \be_j & (y_{ijab} \in 
\Gamma_{\an,r}),
\end{align*}
then $\min_j \{w_r(y_{ijab})\} \geq \min_j\{w_r(x_j)\} + w(b!)$
(i.e., the same bound as for the trivial connection with
$\be_1, \dots, \be_n$ horizontal).
This gives the bound
\begin{align*}
w_r(N_J) &\geq w(J!) + \sum_{i=1}^m \left( -w(j_i!) + 
\lfloor j_i/p^h \rfloor w((p^h)!) +
w((j_i -p^h \lfloor j_i/p^h\rfloor)!) \right) \\
&\geq w(J!) -\sum_{i=1}^m w(p)j_i/(p^h(p-1))
\end{align*}
using the fact that $w(j_i!) = \sum_{g=1}^\infty w(p) \lfloor j_i/p^g \rfloor$.
This yields the claim.
\end{proof}

\begin{lemma} \label{L:bound nabla2}
Let $M$ be an $(F,\nabla)$-module over $\Gancon$ or over $\Gcon[\fp]$,
and let $\be_1, \dots, \be_n$ be a basis of $M$.
For each nonnegative integer $g$,
define the $n \times n$ matrices $N_{g,1}, \dots, N_{g,m}$ by $\Delta_i
(F^g \be_l) = \sum_j (N_{g,i})_{jl} (F^g \be_j)$. Then there 
exist $r_1 \in (0,r_0)$ and $c>0$
such that for each nonnegative integer $g$ and for each of $i=1, \dots, m$,
$N_{g,i}$ has entries in $\Gamma_{\an,r_1q^{-g}}$ and
\[
w_{rq^{-g}}(N_{g,i}) \geq g-c \qquad (r \in [r_1/q,r_1]).
\]
Moreover, if $M$ is defined over $\Gcon[\fp]$, we can also ensure that
$w(N_{g,i}) \geq g-c$.
\end{lemma}
\begin{proof}
Define $a_{hi} \in \Gcon$ by the formula
\[
\del_i (x^\sigma) = \sum_{h=1}^m a_{hi} (\del_h x)^\sigma \qquad (x \in \Gcon);
\]
then $w(a_{hi}) \geq 1$ as in
Remark~\ref{R:omega-sigma}. In particular, we can
choose $r_1 \in (0,r_0)$ as in Proposition~\ref{P:naive compare}
such that for $i=1, \dots, m$,
$a_i \in \Gamma_{r_1}$, $w_{r_1}(a_{hi}) \geq 1$, and $N_{0,i}$ has entries in
$\Gamma_{\an,r_1}$. Then the formula
\[
N_{g+1,i} = \sum_{h=1}^m a_{hi} N_{g,h}^\sigma
\]
yields the claim for any $c$ with $\min_i \{\min_{r \in [r_1/q,r_1]} 
\{w_r(N_{0,i})\} \}\geq -c$
and (in case $M$ is defined over $\Gcon[\fp]$)
$\min_i \{w(N_{0,i})\} \geq -c$.
\end{proof}

\begin{lemma} \label{L:convergence2}
Let $M$ be an $(F,\nabla)$-module over $\Gancon$ (resp.\ over $\Gcon[\fp]$).
Then for any $u_1, \dots, u_m \in \Gcon$ with $w(u_i) > 0$ for $i=1,
\dots, m$, and any $\bv \in M$,
the series
\[
f(\bv) = \sum_{j_1, \dots, j_m =0}^\infty 
\frac{1}{J!} U^J \Delta^{\underline{J}}(\bv)
\]
converges for the natural topology of $M$,
and the map $\bv \mapsto f(\bv)$ is semilinear for the map
defined by Lemma~\ref{L:convergence1}.
\end{lemma}
\begin{proof}
Pick a basis $\be_1, \dots, \be_n$ of $M$; 
for each nonnegative integer $g$,
define the $n \times n$ matrices $N_{g,1}, \dots, N_{g,m}$ by $\Delta_i
(F^g \be_l) = \sum_j (N_{g,i})_{jl} (F^g \be_j)$. 
By Lemma~\ref{L:bound nabla2}, we can choose $r_1 \in (0,r_0)$ such that
for some $c>0$, $w_{r q^{-g}}(N_{g,i}) \geq g-c$ for all nonnegative
integers $g$ and all $r \in [r_1/q,r_1]$; moreover,
if we are working over $\Gcon[\fp]$, we can
ensure that $w(N_{g,i}) \geq g-c$.

Now choose a positive integer $h$ with
$w(p)/(p^h(p-1)) < 1/2$. Then
by the previous paragraph, for each sufficiently small $r>0$, there exists
a basis $\bv_1, \dots, \bv_n$ of $M$ (depending on $r$)
on which each $\Delta_i$ acts via a matrix $N_i$
with $w_r(N_i) > w((p^h)!)$. By Lemma~\ref{L:bound nabla},
the matrix $N_J$ defined by
\[
\Delta^{\underline{J}}( \bv_l) = \sum_j (N_J)_{jl} \bv_j
\]
satisfies $w_r(N_J) \geq w(J!)-(j_1+\cdots+j_m)/2$.

On the other hand, since $w(u_i) \geq 1$ for $i=1,\dots, m$, we have
that $w_r(u_i) > 1/2$ for $r$ sufficiently small.
We conclude that for each sufficiently small $r>0$, there exists
a basis $\bv_1, \dots, \bv_n$ such that the series defining each
of $f(\bv_1), \dots, f(\bv_n)$ converges under $w_r$.
By Lemma~\ref{L:Leibniz rule} 
and Lemma~\ref{L:convergence1}, for each $\Gamma_{\an,r}$-linear
combination $\bv$ of $\bv_1, \dots, \bv_n$, the series defining
$f(\bv)$ converges under $w_r$. 
By the same token, in case $M$ is defined over $\Gcon[\fp]$, for 
each $\Gamma_{r}[\fp]$-linear
combination $\bv$ of $\bv_1, \dots, \bv_n$, the series
defining $f(\bv)$ converges under $w$.
This yields the desired convergence
of $f$; again, the semilinearity follows from Lemma~\ref{L:Leibniz rule}
and Lemma~\ref{L:convergence1}.
\end{proof}

\begin{prop} \label{P:change Frob}
Let $\sigma_1$ and $\sigma_2$ be Frobenius lifts on $\Gamma$ 
such that $\Gamma$ has enough units with 
respect to each of $\sigma_1$ and $\sigma_2$, and 
for each $z \in L$ nonzero, $\{z\}$ is a unit in $\Gcon$
under both definitions.
(By Proposition~\ref{P:naive compare}, it is equivalent to require
that the definitions of $\Gcon$ with respect to $\sigma_1$ and to $\sigma_2$
coincide.)
Then there is a canonical equivalence of categories between 
$(F, \nabla)$-modules over $\Gancon$ (resp.\ over $\Gcon[\fp]$)
relative to $\sigma_1$
and relative to $\sigma_2$, acting as the identity on the underlying
$\nabla$-modules.
\end{prop}
\begin{proof}
Put  $u_i = \{z_i\}^{\sigma_2}/\{z_i\}^{\sigma_1}-1$.
Let $M$ be a $\nabla$-module admitting a compatible
Frobenius structure $F_1$ relative to $\sigma_1$.
For $\bv \in M$, define
\[
F_2(\bv) = \sum_{j_1, \dots, j_m =0}^\infty
\frac{1}{J!} U^J F_1(\Delta^{\underline{J}}(\bv));
\]
this series converges thanks to Lemma~\ref{L:convergence2}.
Moreover, the result is $\sigma_2$-linear thanks to 
Lemma~\ref{L:Leibniz rule}.
\end{proof}

\begin{remark}
By tweaking the proof of 
Proposition~\ref{P:change Frob}, one can also obtain the
analogous independence from the choice of $\sigma$ for
the category of $(F, \nabla)$-modules over $\Gcon$. We will not
use this result explicitly, though a related construction will occur
in Subsection~\ref{subsec:galois rep}.
\end{remark}

\section{Unit-root $(F, \nabla)$-modules (after Tsuzuki)}
\label{sec:unit-root}

In this section, we give the generalization to fake annuli
of Tsuzuki's unit-root local
monodromy theorem \cite{tsuzuki-unitroot},
variant proofs of which are given by Christol
\cite{christol} and in the
author's unpublished dissertation \cite{me-thesis}.
Our argument here draws on elements of all of these; its specialization to
the case of true annuli constitutes a 
novel (if only slightly so) exposition of Tsuzuki's original result.

\setcounter{theorem}{0}
\begin{convention} \label{conv:finite fake}
Throughout this section, let $E$ denote a nearly monomial field over $k$,
viewed in a fixed fashion as a finite separable
extension of a monomial field over $k$.
We assume that any Frobenius lift $\sigma$ considered
on $\Gamma = \Gamma^E$ is chosen so that
$\Gamma$ has enough units.
In particular,
$\Gamma = \Gamma^E$ and $\Gcon = \Gcon^E$ are nearly admissible
in the sense of Definition~\ref{D:admissible}.
\end{convention}

\subsection{Unit-root $F$-modules}

\begin{defn}
We say an $F$-module $M$ over $\Gamma^E$ or $\Gcon^E$,
with respect to some Frobenius lift $\sigma$, 
is \emph{unit-root} (or \emph{\'etale})
if the map $F: \sigma^* M \to M$ is an isomorphism
(not just an isogeny). 
We say an $(F, \nabla)$-module over $\Gamma^E$
or $\Gcon$ is unit-root if its underlying $F$-module is unit-root.
\end{defn}

We will frequently calculate on such modules in terms of bases, so it is worth
making the relevant equations explicit.
\begin{remark} \label{R:equations}
Assume that $E$ is a monomial field, and fix a coordinate
system for $\Gamma$.
Let $M$ be a $\nabla$-module over $\Gamma$ or $\Gcon$
with basis $\be_1, \dots, \be_n$.
Given $\mu \in L^\dual$, define the $n \times n$ matrix $N_\mu$ by
$\Delta_\mu (\be_l) = \sum_j (N_\mu)_{jl} \be_j$; if we identify
$\bv = c_1 \be_1 + \cdots + c_n \be_n \in M$ with the column
vector with entries $c_1, \dots, c_n$, then we have
\[
\Delta_\mu (\bv) = N_\mu \bv + \del_\mu (\bv).
\]
Given an $F$-module with the same basis $\be_1, \dots, \be_n$, 
define the $n \times n$ matrix $A$
by $F(\be_l) = \sum_j A_{jl} \be_j$; then with the same identification
of $\bv$ with a column vector, we have
\[
F(\bv) = A \bv^\sigma.
\]
In case the Frobenius lift $\sigma$ is standard,
the compatibility of Frobenius and connection structures is
equivalent to the equations
\[
N_\mu A + \del_\mu (A) = q A N_\mu^\sigma \qquad (\mu \in L^\dual);
\]
of course it is only necessary to check this on a basis of $L^\dual$.
\end{remark}

\begin{remark} \label{R:basis change}
It is also worth writing out how the equations in
Remark~\ref{R:equations} transform under change of basis.
First, if $U$ is an invertible $n \times n$ matrix, then
\[
N_\mu A + \del_\mu (A) = 0 \qquad \Longleftrightarrow \qquad
(U^{-1}N_\mu U + U^{-1} \del_\mu (U))(U^{-1} A) + \del_\mu (U^{-1} A) = 0.
\]
Second, in case $\sigma$ is standard, the equations
\[
N_\mu A + \del_\mu (A) = qAN_\mu^\sigma \qquad \mbox{and} \qquad
N'_\mu A' + \del_\mu (A') = q A'(N'_\mu)^\sigma
\]
are equivalent for
\begin{align*}
N'_\mu &= U^{-1} N_\mu U + U^{-1} \del_\mu (U) \\
A' &= U^{-1} A U^\sigma.
\end{align*}
\end{remark}

\subsection{Unit-root $F$-modules and Galois representations}
\label{subsec:galois rep}

We now consider unit-root $F$-modules over $\Gamma$, obtaining the
usual Fontaine-style setup. 

\begin{lemma} \label{L:solve eq}
Let $\ell$ be a separably closed field of characteristic $p>0$,
and let $\tau$ denote the $q$-power Frobenius on $\ell$.
Let $A$ be an invertible $n \times n$ matrix over $\ell$. 
\begin{enumerate}
\item[(a)]
There exists an invertible $n \times n$ matrix $U$ over $\ell$ such that
$U^{-1} A U^\tau$ is the identity matrix.
\item[(b)]
For any $1 \times n$ column vector $\bv$ over $\ell$, 
there are exactly $q^n$ distinct $1 \times n$ column vectors $\bw$ over
$\ell$ for which $A\bw^\tau - \bw = \bv$.
\end{enumerate}
\end{lemma}
\begin{proof}
Part (a) is \cite[Proposition~1.1]{katz}; part (b) is an easy corollary of (a).
\end{proof}

We next introduce 
a ``big ring'' over which unit-root $F$-modules over $\Gamma$ 
can be trivialized.
\begin{defn}
Let $\tilde{\Gamma}$ be the $\pi$-adic
completion of the maximal unramified extension
of $\Gamma$; then any Frobenius lift $\sigma$ on $\Gamma$ extends uniquely
to $\tilde{\Gamma}$, and the derivation $d$ extends uniquely to a derivation
$d: \tilde{\Gamma} \to (\Omega^1_{\Gamma/\gotho} \otimes_{\Gamma} \tilde{\Gamma})$.
Likewise, any $\nabla$-module $M$ over $\Gamma$ induces a connection
$\nabla: (M \otimes_\Gamma \tilde{\Gamma}) \to (M \otimes_\Gamma 
\Omega^1_{\Gamma/\gotho} \otimes_{\Gamma} \tilde{\Gamma})$.
%Let $W$ be the $\pi$-adic completion of the direct limit
%$\tilde{\Gamma} \stackrel{\sigma}{\to} \tilde{\Gamma} 
%\stackrel{\sigma}{\to} \cdots$.
Let $\tilde{\gotho}_q$ be the fixed subring of $\tilde{\Gamma}$
under $\sigma$; this is a complete discrete valuation ring with residue field
$\FF_q$ and maximal ideal generated by $\pi$.
\end{defn}

\begin{prop} \label{P:invariant basis}
Let $M$ be a unit-root $F$-module over $\tilde{\Gamma}$. Then $M$
admits an $F$-invariant basis.
\end{prop}
\begin{proof}
Applying Lemma~\ref{L:solve eq}(a) produces a basis which is fixed modulo $\pi$.
Given a basis fixed modulo $\pi^n$, correcting it to a basis fixed modulo
$\pi^{n+1}$ amounts to solving a set of vector equations of the form
of Lemma~\ref{L:solve eq}(b). The resulting sequence of bases converges
to the desired $F$-invariant basis.
\end{proof}

\begin{defn} \label{D:equiv}
Assume that $\FF_q \subseteq k$.
Given a unit-root $F$-module
$M$ over $\Gamma$, let $D_\Gamma(M)$ denote the set of $F$-invariant
elements of $M \otimes_{\Gamma} \tilde{\Gamma}$;
then $D_\Gamma(M)$ is a finite free $\tilde{\gotho}_q$-module equipped 
with a continuous action of $G = \Gal(E^{\sep}/E)$.
By Proposition~\ref{P:invariant basis}, the natural map
$D_{\Gamma}(M) \otimes_{\tilde{\gotho}_q} \tilde{\Gamma} 
\to M \otimes_\Gamma \tilde{\Gamma}$
is an isomorphism.
Conversely, given a finite free 
$\tilde{\gotho}_q$-module $N$ equipped with a
continuous action of $G$,
let $V(N)$ denote the set of $G$-invariant elements of $N \otimes_{\tilde{\gotho}_q} 
\tilde{\Gamma}$;
by Galois descent, the natural map
$V(N) \otimes_\Gamma \tilde{\Gamma} \to N 
\otimes_{\tilde{\gotho}_q} \tilde{\Gamma}$
is an isomorphism. The functors $D_\Gamma$ and $V$ thus exhibit
equivalences of categories between unit-root $F$-modules over $\Gamma$
and finite free $\tilde{\gotho}_q$-modules equipped with continuous $G$-action.
\end{defn}

So far in this subsection,
we have considered only unit-root $F$-modules over $\tilde{\Gamma}$,
without connection structure. The reason is that the connection does not
really add any extra structure in this case.
\begin{prop} \label{P:unique}
Let $M$ be a unit-root $F$-module over $\Gamma$ (resp.\ over 
$\tilde{\Gamma}$).
Then there is a unique integrable connection $\nabla: M \to M \otimes 
\Omega^1$ compatible with $F$.
\end{prop}
\begin{proof}
We first check existence and uniqueness for the connection, without
worrying about integrability. Let $\nabla_0: M \to M \otimes \Omega^1$
be any connection (not necessarily integrable). 
Then a map $\nabla: M \to M \otimes \Omega^1$ is a connection
if and only if $\nabla - \nabla_0$ is a 
$\Gamma$-linear map from $M$ to $M \otimes \Omega^1$,
i.e., if it corresponds to an element of $M^\dual \otimes M \otimes \Omega^1$.
Moreover, $\nabla$ is compatible with $F$ if and only if
\[
(\nabla -\nabla_0) F - (F \otimes d\sigma)(\nabla-\nabla_0) = 
(F \otimes d\sigma) \nabla_0 - \nabla_0 F;
\]
in other words, we can write down a particular $\bw \in M^\dual
\otimes M \otimes \Omega^1$ such that $\nabla$ is $F$-equivariant
if and only if $\nabla - \nabla_0$ corresponds to an element $\bv \in M^\dual
\otimes M \otimes \Omega^1$ with $\bv - F\bv = \bw$. By
Remark~\ref{R:omega-sigma},
$M^\dual \otimes M \otimes \Omega^1$ can be written as 
a twist $N(1)$ for some $F$-module $N$ over $\Gamma$ (resp.\
over $\tilde{\Gamma}$); hence the
series $\bw + F\bw + F^2 \bw + \cdots$ converges in $M^\dual \otimes M
\otimes \Omega^1$ to the unique solution of $\bv - F\bv = \bw$.

It remains to prove that the unique connection $\nabla$ compatible with
$F$ is in fact integrable. It is enough to check this over $\tilde{\Gamma}$;
moreover, it is enough to exhibit a single integrable connection compatible
with $F$, as this must then coincide with the connection constructed above.
To do this, we apply Proposition~\ref{P:invariant basis} to produce an
$F$-invariant basis $\be_1, \dots, \be_n$ of $M$, then set
\[
\nabla(c_1 \be_1 + \cdots + c_n \be_n) = \be_1 \otimes dc_1 + \cdots +
\be_n \otimes dc_n;
\]
this map is easily seen to be an integrable connection compatible
with $F$.
\end{proof}

\begin{remark} \label{R:change Frob}
For true annuli, the construction of Definition~\ref{D:equiv}
is due to Fontaine \cite[1.2]{fontaine}.
In general, one consequence of the construction
is that the categories of unit-root $F$-modules over $\Gamma$
relative to two different Frobenius lifts are canonically equivalent,
since the category of finite free $\tilde{\gotho}_q$-modules equipped with 
continuous $G$-action does not depend on the choice of the Frobenius lift.
In fact, one may even change the choice of the underlying
Frobenius lift $\sigma_K$, as long as it does not change what 
$\tilde{\gotho}_q$ is.
Note that the same is true of unit-root $(F, \nabla)$-modules; more
precisely, if a $\nabla$-module $M$ over $\Gamma$
admits a unit-root Frobenius structure
for one Frobenius lift, it admits a unit-root Frobenius structure
for any Frobenius lift. That is because the connection on $M$
can be recovered from $D_\Gamma(M)$, by specifying that elements of 
$D_\Gamma(M)$ are horizontal.
\end{remark}

\subsection{Positioning Frobenius}

It will be useful to prove a positioning lemma for 
elements of $k((L))_\lambda$.

\begin{lemma} \label{L:positioning}
Assume that the field $k$ is algebraically closed.
Suppose $x \in k((L))_\lambda$ cannot be written as $a^q - a$ for any
$a \in k((L))_\lambda$. Then there exists $c>0$ such that
for any nonnegative integer $i$ and
any $y \in k((L))_\lambda$ with $x - y + y^q \in k((q^i L))_\lambda$,
we have $v_\lambda(x - y + y^q) \leq -cq^i$.
\end{lemma}
\begin{proof}
Clearly there is no harm in replacing $x$ by $x - y_0 + y_0^q$
for any $y_0 \in k((L))_\lambda$.
In particular, write $x = \sum_{z \in L} c_z \{z\}$, and let
$x_-, x_0, x_+$ be the sum of $c_z \{z\}$ over those $z$ with
$\lambda(z)$ negative, zero, positive, respectively. Since $k$
is algebraically closed, we have $x_0 = y - y^q$ for some $y \in k$.
Since $v_\lambda(x_+) > 0$, we have $x_+ = y - y^q$ for
$y = x_+ + x_+^q + x_+^{q^2} + \cdots$.
We may thus reduce to the case $x = x_-$; in particular,
$x$ has finite support.

For $z \in L$ nonzero,
let $i(z)$ denote the largest integer $i$ such that
$z/q^i \in L$. Set
\[
y_1 = \sum_{z \in L, \lambda(z) < 0} \sum_{i=1}^{i(z)} 
(c_z \{z\})^{q^{-i}},
\]
so that $x_1 = x + y_1 - y_1^q$ is supported on $L \setminus qL$. We cannot
have $x_1 = 0$, or else we could have written $x = a^q - a$ for some
$a \in k((L))_\lambda$. There must thus be a smallest (under $\lambda$)
element $z$ of the support of $x_1$.
For any nonnegative integer $i$
and any $y \in k((L))_\lambda$ with $x - y + y^q \in k((q^i L))_\lambda$,
the support of $x - y + y^q$ must contain $q^{i+j} z$ for some
nonnegative integer $j$, and so 
$v_\lambda(x - y + y^q) \leq -\lambda(z) q^i$, as desired.
\end{proof}

\subsection{Successive decimation}

We now give a version of Tsuzuki's 
construction for solving $p$-adic differential equations.
\begin{convention} \label{conv:coordinate}
Throughout this subsection, assume that $k$ is algebraically closed
and that $E$ is a monomial field over $k$,
and fix a coordinate system on $\Gamma$.
Also assume $\sigma$ is a standard Frobenius lift;
we may then safely confound the na\"\i ve and Frobenius-based
partial valuations.
\end{convention}

\begin{defn}
For $\mu \in L^\dual$ nonzero,
write $L_\mu$ for the sublattice of $z \in L$ for which
$\mu(z) \in p\ZZ$.
\end{defn}

\begin{lemma} \label{L:decimate0}
Suppose that $A$ is an invertible $n \times n$ matrix over
$\Gamma$ with $w(A-I_n) > 0$, 
that $N_\mu$ is an $n \times n$ matrix over
$\Gamma$ supported on $L_\mu$,
and that $N_\mu A + \del_\mu (A) = qAN_\mu^\sigma$. Then $A$ is supported on
$L_\mu$.
\end{lemma}
\begin{proof}
Suppose the contrary; write $A = B + C$ with $B$ supported on $L_\mu$
and $h = w(C)$ as large as possible, so in particular $h > 0$. 
Write $C = \sum_{z \in L} C_z \{z\}$.
Since $h$ is as large as possible, there exists $z \in L \setminus L_\mu$
such that $w(C_z) = h$; the coefficient of $\{z\}$ in $\del_\mu (C)$ then also
has valuation $h$. 
However, in the equality
\[
\del_\mu (C) = (qBN_\mu^\sigma - N_\mu B - \del_\mu (B)) + (q C N_\mu^\sigma
- N_\mu C),
\]
the first term in parentheses is supported on $L_\mu$ while the second
term has valuation strictly greater than $h$
(since $w(A-I_n) > 0$ forces $w(N_\mu) > 0$).
This contradiction yields the desired result.
\end{proof}

\begin{lemma} \label{L:decimate1}
Pick $\mu \in L^\dual$ nonzero.
Given $r>0$, let $N_\mu$ be an $n \times n$ matrix over
$\Gamma_r$ such that $w(N_\mu) > 0$ and 
$w_r(N_\mu) > 0$.
Then for any $s \in (0,r)$,
there exists an invertible $n \times n$ matrix $U$ over $\Gamma_s$
such that $w_s(U-I_n) > 0$, $w(U-I_n) \geq w(N_\mu)$,
and $U^{-1} N_\mu U + U^{-1} \del_{\mu} (U)$ is supported on $L_\mu$.
\end{lemma}
\begin{proof}
Define a sequence $U_0, U_1, \dots$
of invertible matrices over $\Gamma_{r}$ with $w(U_j - I_n) \geq w(N_\mu)$
and $w_r(U_j - I_n) \geq w_r(N_\mu)$, as follows. 
Start with $U_0 = I_n$. 
Given $U_j$, put $N_{j} = U_j^{-1} N_\mu U_j + U_j^{-1} \del_\mu (U_j)$.
Write $N_{j} = \sum_{z\in L} N_{j,z} \{z\}$, let $X_{j}$ be the sum of $\mu(z)^{-1}
N_{j,z} \{z\}$ over all $z \in L \setminus L_\mu$,
and put $U_{j+1} = U_j(I_n - X_j)$.

For $j>0$, if $w(U_j - I_n) < \infty$, one sees that 
$w(U_{j+1} - I_n) > w(U_j - I_n)$.
Hence the $U_j$ converge $\pi$-adically; 
since they all satisfy $w_r(U_j-I_n) 
\geq w_r(N_\mu) > 0$, the $U_j$ converge under $w_s$ to a limit $U$ satisfying
$w_s(U-I_n) \geq w_s(N_\mu)$. In particular, $U$ is invertible and
$U^{-1}N_\mu U + U^{-1} \del_\mu (U)$ is supported on $L_\mu$.
(Compare \cite[Lemma~6.1.4]{tsuzuki-unitroot}
and \cite[Lemma~5.1.3]{me-thesis}.)
\end{proof}

\begin{lemma} \label{L:decimate3}
Let $N_1, \dots, N_m$ be $n \times n$ matrices over
$\Gamma_r$, such that $w(N_i) > 0$ and $w_r(N_i) > 0$
for $i=1, \dots, m$. Suppose that there exists an invertible matrix
$A$ over $\Gamma$ with $w(A-I_n) > 0$, 
such that 
\[
N_i A + \del_i (A) = qA N_i^\sigma \qquad (i=1, \dots, m).
\]
Then 
for any $s \in (0,r)$,
there exists an invertible matrix $U$ over $\Gamma_s$,
such that $w(U-I_n) \geq \min_i \{w(N_i)\}$,
$w_s(U-I_n) \geq \min_i \{w_s(N_i)\}$,
and $U^{-1} AU^\sigma$ is supported on $pL$.
\end{lemma}
\begin{proof}
For $i=0, \dots, m$, let $S_i$ be the sublattice of $z \in L$
such that $\mu_j(z) \in p\ZZ$ for $j = 1,\dots, i$.
Pick $s_1,\dots,s_{m-1}$ with $0 < s < s_{m-1} < \cdots < s_1 < r$,
and put $s_0 = r$ and $s_{m} = s$.
Put $U_0 = I_n$. Given 
$U_i$ invertible over $\Gamma_{s_i}$
such that $A_i = U_i^{-1} A U_i^\sigma$ is supported on $S_i$,
note that $M_i = U_i^{-1} N_i U_i + U_i^{-1} \del_i (U_i)$ 
satisfies the equation $M_i A_i + \del_i (A_i) = q A_i M_i^\sigma$.
We may then argue (as in the proof of Proposition~\ref{P:unique})
that $M_i$ is congruent to a matrix supported on $S_i$ modulo
successively larger powers of $\pi$.

Since $M_i$ is supported on $S_i$, we may
 apply Lemma~\ref{L:decimate1} to produce $U_{i+1} = U_i V_i$, with
$V_i$ supported on $S_i$, 
such that $U_{i+1}^{-1} N_i U_{i+1} + U_{i+1}^{-1} \del_i (U_{i+1})$
is supported on $S_{i+1}$. Then
$U_{i+1}^{-1} A U_{i+1}^\sigma = V_i^{-1} A_i V_i^\sigma$ 
is supported on $S_{i}$;
by Lemma~\ref{L:decimate0},
$U_{i+1}^{-1} A U_{i+1}^\sigma$ is also supported on $S_{i+1}$. 
Thus the iteration continues, and we may set $U = U_{m}$.
\end{proof}

We are now ready for the decisive step, analogous to
\cite[Lemma~5.2.4]{tsuzuki-unitroot}.
\begin{prop} \label{P:tsuzuki}
Let $N_1, \dots, N_m$ be  $n \times n$ 
matrices over $\Gcon$ with $w(N_i) > 0$ for $i=1,\dots,m$.
Suppose that there exists an invertible $n \times n$ matrix
$A$ over $\Gamma$ such that $w(A-I_n) > w(p)/(p-1)$ and
\[
N_i A + \del_i (A) = q AN_i^\sigma \qquad (i=1, \dots, m).
\]
Then there exists an invertible $n \times n$ matrix $U$ over $\Gamma$
such that $AU^\sigma = U$.
\end{prop}
\begin{proof}
Suppose the contrary; then there exists some smallest integer 
$h > w(p)/(p-1)$ such that the equation $U^{-1} A U^\sigma \equiv
I_n \pmod{\pi^{h+1}}$ cannot be solved for $U$ invertible over
$\Gamma$. Since $\Gcon$ is $\pi$-adically dense in $\Gamma$, we may
change basis over $\Gcon$ to reduce to the case where
 $h = w(A - I_n)$
and we cannot write the reduction of $\pi^{-h}(A-I_n)$ modulo $\pi$
in the form $B - B^\sigma$.

Choose $r_0 >0$ such that for $i=1, \dots, m$, $N_i$ has entries in 
$\Gamma_{r_0}$ and $w_{r_0}(N_i) > 0$. Since $h > w(p)/(p-1)$,
we have $hp/(h+w(p)) > 1$; we can thus choose $c$ with 
$1 < c < hp/(h+w(p))$. Write $r_j = r_0 p^{-j} c^{j}$.

Define a sequence $U_0, U_1, \dots$ of invertible matrices over $\Gcon$
as follows. Start with $U_0 = I_n$. For $j \geq 0$, suppose that
we have constructed an invertible matrix $U_j$ over $\Gamma_{r_j}$
with the following properties:
\begin{enumerate}
\item[(a)]
$w(U_j - I_n) \geq h$ and $w_{r_j}(U_j - I_n) > 0$;
\item[(b)]
$A_j = U_j^{-1} A U_j^\sigma$ 
and $N_{i,j} = U_j^{-1} N_i U_j + U_j^{-1} \del_i (U_j)$
are supported on $p^j L$ for $i=1, \dots, m$;
\item[(c)]
$w_{r_j}(N_{i,j}) > j w(p)$.
\end{enumerate}
Write $A'_j$ and $N'_{i,j}$ for the matrices $A_j$ and $p^{-j} N_{i,j}$
viewed in $\Gamma^{E_j}$ for $E_j = k((p^j L))_\lambda$,
put $\mu'_i = p^{-j} \mu_i$, and let
$\del'_i$ be the derivation on $\Gamma^{E_j}$ corresponding to
$\mu'_i$. Then $w_{r_j}(N'_{i,j}) > 0$, and
\[
N'_{i,j} A'_j + \del'_i (A'_j) = q A'_j (N'_{i,j})^\sigma.
\]
Put $s = r_j (h+w(p))c/(hp) < r_j$,
and apply Lemma~\ref{L:decimate3} 
to produce $U_{j+1}$ over $\Gamma_s$ supported on
$p^j L$, with
$w(U_{j+1} - I_n) \geq h$ and $w_s(U_{j+1} - I_n) > 0$,
such that $A_{j+1} = U_{j+1}^{-1} A_j U_{j+1}^\sigma$ is supported
on $p^{j+1} L$.

Since $r_{j+1} = sh/(h + w(p)) < s$, (a) is satisfied again.
From the equation 
\begin{equation} \label{eq:decimate}
N_{i,j+1} A_{j+1} + \del_i (A_{j+1}) = q A_{j+1}
N_{i,j+1}^\sigma
\end{equation}
(a consequence of Remark~\ref{R:basis change}), we see that
each $N_{i,j+1}$ is also supported on $p^{j+1}L$
(the argument is as in the proof of Lemma~\ref{L:decimate3}).
Hence (b) is satisfied again.

To check (c), note that on one hand,
\eqref{eq:decimate} and the fact that $w(\del_i (A_{j+1})) \geq h
+ (j+1) w(p)$ imply that $w(N_{i,j+1}) \geq h + (j+1)w(p)$.
On the other hand, the facts that $w_{s}(N_{i,j}) > w_{r_j}(N_{i,j}) > jw(p)$,
$U_{j+1}$ is supported on $p^j L$, and
$w_s(U_{j+1} - I_n) > 0$ imply that $w_s(N_{i,j+1}) > jw(p)$,
and so
\begin{align*}
w_{r_{j+1}}(N_{i,j+1}) 
&= \min_{m \geq h+(j+1)w(p)} \{r_{j+1} v_m(N_{i,j+1}) + m\} \\
&\geq \frac{r_{j+1}}{s} w_s(N_{i,j+1}) + (h+(j+1)w(p)) \left( 1 - 
\frac{r_{j+1}}{s} \right) \\
&> \frac{r_{j+1}}{s} (jw(p)) + (h+(j+1)w(p)) \left( 1 - 
\frac{r_{j+1}}{s} \right) \\
&= \frac{h}{h+w(p)} (jw(p)) + (h+(j+1)w(p)) \frac{w(p)}{h + w(p)} \\
&= (j+1) w(p).
\end{align*}
Hence (c) is satisfied again,
and the iteration may continue.

Note that if we take $X = U_j - I_n$, then $A - X + X^\sigma \equiv A_j
\pmod{\pi^h}$. This means that on one hand, $A - X + X^\sigma$
is congruent modulo $\pi^h$ to a matrix supported on $p^j L$,
and on the other hand,
\begin{align*}
v_h(A - X + X^\sigma) &\geq \min\{v_h(A), v_h(X), v_h(X^\sigma)\} \\
&\geq \min\{v_h(A), -h r_0^{-1} p^{j+1} c^{-j} \}
\end{align*}
since $w_{r_j}(X) > 0$. However, since $c>1$, this last
inequality contradicts Lemma~\ref{L:positioning} for $j$ large.
This contradiction means that our original assumption was
incorrect, i.e., the desired matrix $U$ does exist, as desired.
\end{proof}

\begin{remark}
In his setting, Tsuzuki actually proves a stronger result 
\cite[Proposition~6.1.10]{tsuzuki-unitroot} that produces solutions
of $p$-adic differential equations even without a Frobenius structure.
As noted in Remark~\ref{R:cannot integrate}, one cannot hope to
do likewise in our setting.
\end{remark}

\subsection{Trivialization}

We now begin reaping the fruits of our labors, first in a 
restricted setting. Note that Convention~\ref{conv:coordinate} is 
no longer in force.
\begin{prop} \label{P:triv Frob}
Suppose that the field $k$ is algebraically closed.
Let $M$ be a $\nabla$-module over $\Gcon$ which becomes
an $(F, \nabla)$-module over $\Gamma$. 
Then as a representation of $G = \Gal(E^{\sep}/E)$,
$D_\Gamma(M)$ has finite image; moreover, if
$D_\Gamma(M)$ is trivial modulo $\pi^m$ for some
integer $m > w(p)/(p-1)$, then $D_\Gamma(M)$ is trivial.
\end{prop}
\begin{proof}
We treat the second assertion first. Suppose that
$D_\Gamma(M)$ is trivial modulo $\pi^m$ for some
integer $m > w(p)/(p-1)$.
By Remark~\ref{R:change Frob}, we can change the choice of the Frobenius
lift without affecting the fact that $\nabla$ admits a compatible
Frobenius structure over $\Gamma$, or that $D_\Gamma(M)$ is trivial
modulo $\pi^m$. In particular, we may assume that
$\sigma$ is a standard Frobenius lift; we may then choose a coordinate
system for $\Gamma$
to drop back into the purview of Convention~\ref{conv:coordinate}.

Given a basis $\be_1, \dots, \be_n$ of $M$, define
the $n \times n$ matrices $A, N_1, \dots, N_m$ by
\begin{align*}
F(\be_l) &= \sum_j A_{jl} \be_j \\
\Delta_i (\be_l) &= \sum_j (N_i)_{jl} \be_j;
\end{align*}
as in Remark~\ref{R:equations}, we then have
$N_i A + \del_i (A) = q A N_i^\sigma$ for $i=1, \dots, m$.
By hypothesis, we can arrange to have $w(A-I_n) > w(p)/(p-1)$;
we may thus apply Proposition~\ref{P:tsuzuki} to produce an invertible
$n \times n$ matrix $U$ over $\Gamma$ with $A U^\sigma = U$.
Writing $\bv_l = \sum_j U_{jl} \be_j$, we then have $F\bv_j = \bv_j$
for $j=1, \dots, n$; that is, $M \otimes_{\Gcon} \Gamma$ 
admits an $F$-invariant basis, so $D_\Gamma(M)$ is trivial.

We now proceed to the first assertion.
Pick an integer $m$ with $m > w(p)/(p-1)$.
Let $E'$ be the fixed field of the kernel of the action of
$G$ on $D_\Gamma(M)/\pi^m D_\Gamma(M)$.
Then by what we have just shown, the restriction of 
$D_\Gamma(M)$ to $\Gal(E^{\sep}/E')$ is trivial; hence
$D_\Gamma(M)$, as a representation of $G$, has finite image.
\end{proof}

Finally, we give the analogue of Tsuzuki's unit-root monodromy
theorem \cite[Theorems~4.2.6 and~5.1.1]{tsuzuki-unitroot}.
\begin{theorem} \label{T:tsuzuki}
Let $M$ be a unit-root $(F,\nabla)$-module over $\Gcon$. Then
there exists a finite separable extension $E'$ of $E$ such that
$M \otimes_{\Gcon} \Gcon^{E'}$ admits a basis of elements which are horizontal,
and also $F$-invariant in case $k$ is algebraically closed.
\end{theorem}
\begin{proof}
Suppose $k$ is algebraically closed; by
Proposition~\ref{P:triv Frob}, 
for some finite separable extension $E'$ of $E$,
there exists a basis of $M \otimes_{\Gcon} \Gamma'$ (for $\Gamma' = \Gamma^{E'}$) 
consisting of
horizontal $F$-invariant elements. Since $M$ is unit-root, we may
apply \cite[Lemma~5.4.1]{me-slope} to deduce that
any $F$-invariant element of $M \otimes_{\Gcon} \Gamma'$ actually belongs
to $M \otimes_{\Gcon} \Gcon'$; this yields the claim.

For $k$ general, let $E'$ denote the completion of the compositum
of $E$ and $k^{\alg}$ over $k$. Then the restriction of
$D_{\Gamma}(M)$ to $\Gal((E')^{\sep}/E')$ has finite image by
Proposition~\ref{P:triv Frob}. By a standard approximation argument,
we can replace $E$ by a finite separable 
extension in such a way as to trivialize
the action of the resulting $\Gal((E')^{\sep}/E')$; the
resulting action of $\Gal(k^{\sep}/k)$ is trivial by Hilbert 90.
This yields the claim. (Alternatively, one may proceed as in
\cite[Proposition~6.11]{me-local} to reduce the case of $k$ general
to the case of $k$ algebraically closed.)
\end{proof}

\section{Monodromy of $(F,\nabla)$-modules}
\label{sec:monodromy}

In this section, we recall the slope filtration theorem of
\cite{me-local} (in the form presented in \cite{me-slope}), and combine
it with the unit-root monodromy theorem (Theorem~\ref{T:tsuzuki})
to obtain the $p$-adic local monodromy theorem for fake annuli
(Theorem~\ref{T:local mono}).

Throughout this section,
we retain Convention~\ref{conv:finite fake}.

\subsection{Isoclinicity}

\begin{defn}
Let $M$ be an $F$-module of rank 1 over $\Gancon$, and let $\bv$ be a
generator of $M$. Then $F\bv = r\bv$ for some $r \in \Gancon$ which is 
a unit, that is, $r \in \Gcon[\fp]$. In particular, $w(r)$ is well-defined;
it also does not depend on $r$, since changing the choice of generator
multiplies $r$ by $u^\sigma/u$ for some $u \in \Gcon[\fp]$, whereas
$w(u^\sigma/u) = 0$. We call the integer $w(r)$ the \emph{degree} of $M$,
and denote it by $\deg(M)$; if $M$ has rank $n>1$, we define the degree of $M$
as $\deg(\wedge^n M)$. We write $\mu(M) = \deg(M)/\rank(M)$ and call
it the \emph{slope} of $M$.
\end{defn}

\begin{defn}
An $F$-module $M$ over $\Gcon[\fp]$ is \emph{unit-root} (or \emph{\'etale})
if it
contains an $F$-stable $\Gcon$-lattice which forms a unit-root
$F$-module over $\Gcon$. Note that if $M$ is a unit-root 
$(F, \nabla)$-module over $\Gcon[\fp]$, 
then any unit-root $\Gcon$-lattice is stable
under $\nabla$ (as can be seen by applying Frobenius repeatedly).
\end{defn}

\begin{defn} \label{D:isoclinic}
An $F$-module $M$ over $\Gcon[\fp]$ is \emph{isoclinic of slope $s$}
if there exist integers $c$ and $d$ with $c/d = s$
such that $([d]_* M)(-c)$ is unit-root; note that necessarily
$s = \mu(M)$.
An $F$-module $M$ over $\Gancon$ is \emph{isoclinic of slope $s$}
if it is the base extension of an isoclinic $F$-module of slope $s$
over $\Gcon[\fp]$; the base extension from isoclinic $F$-modules
of a given slope over $\Gcon[\fp]$ to isoclinic $F$-modules of that
slope over $\Gancon$ is an equivalence of categories
\cite[Theorem~6.3.3]{me-slope}.
\end{defn}
\begin{remark}
Note that this is not the definition of isoclinicity used in
\cite{me-slope}, but it is equivalent to it thanks to
\cite[Proposition~6.3.5]{me-slope}.
\end{remark}

The base extension functor mentioned above also behaves nicely
with respect to connections; see \cite[Proposition~7.1.7]{me-slope}.
\begin{prop} \label{P:isoclinic nabla}
Let $M$ be an isoclinic $F$-module over $\Gcon[\fp]$.
Suppose that $M \otimes \Gancon$, with its given Frobenius, 
admits the structure of an
$(F, \nabla)$-module. Then $M$, with its given Frobenius,
 already admits the structure of an
$(F, \nabla)$-module.
\end{prop}

\subsection{Slope filtrations and local monodromy}

The slope filtration theorem can be stated as follows;
see Remark~\ref{R:monodromy thm} for a precise citation.
\begin{theorem} \label{T:slope filt}
Let $M$ be an $F$-module over $\Gancon$. Then there exists a unique
filtration $0 = M_0 \subset M_1 \subset \cdots \subset M_l = M$
by saturated $F$-submodules with the following properties.
\begin{enumerate}
\item[(a)]
For $i=1, \dots, l$, the quotient $M_i/M_{i-1}$ is isoclinic of some
slope $s_i$.
\item[(b)]
$s_1 < \cdots < s_l$.
\end{enumerate}
\end{theorem}
\begin{defn}
In Theorem~\ref{T:slope filt}, we refer to the numbers
$s_1, \dots, s_l$ as the \emph{Harder-Narasimhan slopes} (or \emph{HN slopes}
for short) of $M$, viewed as a
multiset in which $s_i$ occurs with multiplicity $\rank(M_i/M_{i-1})$.
See  \cite[\S~4.6]{me-slope} for more on the calculus
of the HN slopes.
\end{defn}

The relevance of the slope filtration theorem
to $(F, \nabla)$-modules comes via the following fact
\cite[Proposition~7.1.6]{me-slope}.
\begin{prop} \label{P:nabla stable}
Let $M$ be an $(F, \nabla)$-module over $\Gancon$. Then each
step of the filtration of Theorem~\ref{T:slope filt} is an
$(F, \nabla)$-submodule.
\end{prop}

Using the slope filtration theorem, we easily obtain the
$p$-adic local monodromy theorem.

\begin{theorem}[$p$-adic local monodromy theorem] \label{T:local mono}
Let $M$ be an $(F, \nabla)$-module over $\Gancon$.
Then $M$ is quasi-unipotent; moreover, if $M$ is isoclinic, then $M$
is quasi-constant.
\end{theorem}
\begin{proof}
Let $0 = M_0 \subset M_1 \subset \cdots \subset M_l = M$
be the filtration of the underlying $F$-module of $M$
given by Theorem~\ref{T:slope filt};
by Proposition~\ref{P:nabla stable}, this is also a filtration
of $(F, \nabla)$-submodules. From the definition of isoclinicity
plus Proposition~\ref{P:isoclinic nabla},
each successive quotient $M_i/M_{i-1}$ can be written as
$N_i \otimes \Gancon$, where $N_i$ is an $(F, \nabla)$-module
over $\Gcon[\fp]$ whose underlying $F$-module is isoclinic.

It suffices to check that each $N_i$ is quasi-constant; since that
condition does not depend on the Frobenius structure, we may check
after applying $[d]_*$ and twisting. We may thus reduce to the
case where $N_i$ is unit-root; in that case,
Theorem~\ref{T:tsuzuki} asserts that indeed $N_i$ is quasi-constant,
as desired.
\end{proof}

\begin{remark} \label{R:monodromy thm}
For true annuli, Theorem~\ref{T:local mono}
is what is normally called the ``$p$-adic (local) monodromy theorem''.
The proof here, restricted to that case, is essentially the same as
in \cite[Theorem~6.12]{me-local},
except that the invocation of the slope filtration theorem
\cite[Theorem~6.10]{me-local} is replaced with the more refined form
\cite[Theorem~6.4.1]{me-slope}.
Proofs in the true annuli case have also been given by
Andr\'e \cite[Th\'eor\`eme~7.1.1]{andre}
and Mebkhout \cite[Corollaire 5.0-23]{mebkhout};
these rely not on a close analysis of Frobenius (as in the slope filtration
theorem), but on the close analysis of connections on annuli given by
the $p$-adic index theorem of Christol-Mebkhout \cite{chris-meb4}.
As per Remark~\ref{R:cannot integrate}, it seems unlikely that
such an approach can be made to work in the fake annuli setting, at least
without integrating Frobenius structures into the analysis.
\end{remark}

\begin{remark} \label{R:standard}
In case $k$ is algebraically closed, one can refine the conclusion
of Theorem~\ref{T:local mono}.
Namely, given an constant 
$(F, \nabla)$-module, the $K$-span of the horizontal sections form an
$F$-module over $K$, to which we may apply the classical Dieudonn\'e-Manin 
theorem; the result is a decomposition of the given
$(F, \nabla)$-module into pieces of the form
$[d]^* \Gancon(c)$. Such pieces
are called \emph{standard} in \cite{me-local} and \cite{me-slope}.
\end{remark}

\section{Complements}
\label{sec:complements}

In this section, 
we gather some consequences of the $p$-adic local monodromy theorem for
fake annuli. These generalize known consequences of the ordinary $p$LMT:
calculation of some extension groups, local duality, and
full faithfulness of overconvergent-to-convergent restriction.

Throughout this section,
retain Convention~\ref{conv:finite fake}.

\subsection{Kernels and cokernels}

We calculate some Hom and Ext groups in the category of
$(F, \nabla)$-modules.

\begin{defn}
For $M$ an $F$-module over some ring, let
$H^0_F(M)$ and $H^1_F(M)$ denote the kernel and cokernel of $F-1$ on
$M$. Note that $\Hom_F(M_1,M_2) = H^0_F(M_1^\dual \otimes M_2)$
and $\Ext^1_F(M_1,M_2) = H^1_F(M_1^\dual \otimes M_2)$.
\end{defn}

\begin{defn}
For $M$ an $(F, \nabla)$-module over some ring,
let $H^0_{F,\nabla}(M)$ be the subgroup of $\bv \in M$ with
$F(\bv) = \bv$ and $\nabla (\bv) = 0$.
Let $H^1_{F, \nabla}(M)$ be the set of pairs $(\bv, \omega)
\in M \times (M \otimes \Omega^1)$ with
\begin{equation} \label{eq:commute1}
\omega + \nabla (\bv) = (F \otimes d\sigma)(\omega), \qquad \nabla_1(\omega) = 0,
\end{equation}
modulo pairs of the form $(F(\bw)-\bw, \nabla (\bw))$ for some $\bw \in M$.
Note that $\Hom_{F,\nabla}(M_1,M_2) = H^0_{F,\nabla}(M_1^\dual \otimes M_2)$
and $\Ext^1_{F,\nabla}(M_1,M_2) = H^1_{F,\nabla}(M_1^\dual \otimes M_2)$.
In particular, by Proposition~\ref{P:change Frob}, we can use any
Frobenius lift $\sigma$ to compute $H^i_{F,\nabla}(M)$.
\end{defn}

\begin{lemma} \label{L:H0}
For $d$ an integer, we have
\[
H^0_{F,\nabla}(\Gancon(d)) = \begin{cases} K_q & d=0 \\ 0 & d \neq 0.
\end{cases}
\]
\end{lemma}
\begin{proof}
Note that $\ker(d: \Gancon \to \Omega^1_{\Gancon/\gotho}) = K$.
Then note that
for $x \in K$ nonzero, $w(x^\sigma \pi^d) = w(x)+d$ can only equal
$w(x)$ for $d=0$.
\end{proof}

\begin{prop} \label{P:H1 positive}
Let $M$ be an $(F, \nabla)$-module over $\Gancon$ whose
HN slopes are all positive. Then $H^1_{F, \nabla}(M) = 0$.
\end{prop}
\begin{proof}
Consider a short exact sequence $0 \to M \to N \to \Gancon \to 0$;
by \cite[Proposition~7.4.4]{me-slope}, the exact sequence splits
if and only if $N$ has smallest HN slope zero. In particular,
this may be checked after enlarging $k$, applying $[a]_*$, and passing from
$k((L))_\lambda$ to a finite separable extension. 
By Theorem~\ref{T:local mono},
this allows
us to reduce to the case where $N$ is a successive extension of
twists of trivial $(F, \nabla)$-modules whose slopes are the HN slopes of $N$. 
If these slopes are all  positive,
then repeated application of
Lemma~\ref{L:H0} implies that the map $N \to \Gancon$ is zero,
which it isn't; hence $N$ has smallest HN slope zero.
It follows that $H^1_{F, \nabla}(M) = 0$, as desired.
\end{proof}

\begin{prop} \label{P:H1}
Assume that $k$ is algebraically closed.
Put $n = w(q)$, let $J$ be the subgroup of $x \in K$ satisfying
$q x^\sigma = \pi^n x$, and
let $z_1, \dots, z_m$ be a basis of $L$.
Then for $d$ an integer, we have
\[
H^1_{F,\nabla}(\Gancon(d)) = 
\begin{cases} J \, \dlog\{z_1\} \oplus \cdots \oplus J \, \dlog\{z_m\} 
& d = -n \\ 0 & d \neq -n. 
\end{cases}
\]
\end{prop}
\begin{proof}
For $d > 0$, Proposition~\ref{P:H1 positive} implies that
$H^1_{F,\nabla}(\Gancon(d)) = 0$; we may thus focus on $d \leq 0$.
First suppose $d=0$. Let $0 \to \Gancon \to M \to \Gancon \to 0$
be a short exact sequence of $(F, \nabla)$-modules.
Then by Theorem~\ref{T:slope filt} and Lemma~\ref{L:H0}, $M$
cannot have any nonzero slopes, so $M$ is isoclinic of slope 0;
as in Definition~\ref{D:isoclinic}, we thus obtain a short exact sequence
$0 \to \Gcon[\fp] \to M_0 \to \Gcon[\fp] \to 0$ of 
$(F, \nabla)$-modules over $\Gcon[\fp]$ from which the original sequence
is obtained by tensoring up to $\Gancon$.
Choose a basis $\bv, \bw$ of $M_0$ such that $\bv$ is an $F$-stable
element of $\Gcon[\fp]$ within $M_0$, $\bw$ maps to an $F$-stable
element under
the map $M_0 \to \Gcon[\fp]$, and $F(\bw) -\bw = c\bv$ with $w(c) > w(p)/(p-1)$.
By Proposition~\ref{P:triv Frob}, for some finite separable extension $E'$
of $E$, $M_0 \otimes \Gamma^{E'}_{\con}[\fp]$ 
admits a basis of $F$-stable elements
$\be_1, \be_2$.
By Lemma~\ref{L:H0}, $\bv$ is a $K_q$-linear combination of $\be_1, \be_2$;
we may thus assume that $\be_1 = \bv$. Similarly, we may assume that
$\be_2$ and $\bw$ have the same image under $M \to \Gcon[\fp]$.
Thus the original exact sequence splits, as desired.

Now suppose $d<0$; we
may assume that $\sigma$ is a standard Frobenius lift.
Write $\Gcon^-$ for the subring of $\Gcon$ of series supported on the set
$\{z \in L: \lambda(z) \leq 0\}$.
We first check that
given a pair $(a, \omega)$ representing an element of 
$H^1_{F,\nabla}(\Gancon(d))$, if
$a \in \Gcon^-[\fp]$ and $a$ 
is supported on $(L \setminus qL) \cup \{0\}$, then
we must have $a \in K$.
Put $\omega = x_1 \dlog\{z_1\} + \cdots + x_m \dlog\{z_m\}$, so that
\[
x_i + \del_i( a) = \pi^d q x_i^\sigma \qquad (i=1, \dots, m).
\]
For $z \in L \setminus qL$, suppose that the coefficient of $\{z\}$
in $a$ is nonzero. Choose $i$ with $\mu_i(z) \neq 0$, so that 
the coefficient of $\{z\}$ in $\del_i (a)$ is nonzero, and
for $j=0,1,\dots$, let $c_j$ be the coefficient of
$\{z\}^{q^j}$ in $x_i$. Since $\sigma$ is standard, the coefficient
of $\{z\}$ in $q x_i^\sigma$ is zero; hence $c_0 \neq 0$.
Moreover, $c_{j+1} = \pi^d q c_j^\sigma$ for $j=0,1,\dots$.
It follows that $w(c_j) = w(c_0) + j(w(q) + d)$ for all $j$; however,
by the definition of $\Gancon$, we must have $\liminf_j (w(c_j)/q^j) > 0$,
contradiction. Thus the coefficient of $\{z\}$ in $a$ is zero for
each $z \in L \setminus qL$;
since $a$ is supported on $(L \setminus qL) \cup \{0\}$, we must
have $a \in K$. 

We next check that if $a \in \Gcon^-[\fp]$, then the pair
$(a, \omega)$ represents the same class in $H^1_{F, \nabla}(\Gancon(d))$
as another pair $(a', \omega')$ with $a' \in \Gcon^-[\fp]$ supported
on $(L \setminus qL) \cup \{0\}$, and hence $a' \in K$ as above.
Write $a = \sum_{z \in L} a_z \{z\}$, and for $j=0,1,\dots$,
write $f_j(a)$ for the sum of
$a_z \{z\}$ over all $z \in q^{j} L \setminus q^{j+1} L$. 
Then the sum
\[
y = \sum_{j=0}^\infty \sum_{l=1}^{j} -(\pi^{\{-l\}})^{d} f_j(a)^{\sigma^{-l}},
\]
where $\pi^{\{0\}} = 1$ and $\pi^{\{l+1\}} = (\pi^{\{l\}})^\sigma \pi$,
converges in $\Gcon^-[\fp]$,
so we can represent the same class in $H^1_{F, \nabla}(\Gancon(d))$
by a pair with first member $a' = a - y + \pi^d y^\sigma$.
Since
\[
a' = a_0 + \sum_{j=0}^\infty (\pi^{\{-j\}})^d f_j(a)^{\sigma^{-j}},
\]
$a'$ is supported on $(L \setminus qL) \cup \{0\}$.

Next, we check that any pair
$(a, \omega)$ represents the same class in $H^1_{F, \nabla}(\Gancon(d))$
as another pair $(a', \omega')$ with $a' \in \Gcon^-[\fp]$.
Write $a = \sum_{z \in L} a_z \{z\}$,
and let $a_+, a_0, a_-$ be the sum of $a_z \{z\}$ over those
$z \in L$ with $\lambda(z)$ positive, zero, negative, respectively.
We can then represent the same class in $H^1_{F, \nabla}(\Gancon(d))$
by a pair with first member $a' = a - y + \pi^d y^\sigma$, for
\[
y = \sum_{i=0}^\infty (\pi^{\{i\}})^d
a_+^{\sigma^i}.
\]
Then $a' = a_0 + a_- \in \Gcon^-[\fp]$.

Combining the previous paragraphs, we find that
every element of $H^1_{F, \nabla}(\Gancon(d))$ is represented by a pair
$(a, \omega)$ with $a \in K$, and consequently
$\pi^{-d} x_i = q x_i^\sigma$. Since $k$ is algebraically closed, we
can force $a=0$; moreover, the resulting class representative is in fact
unique. This yields the desired result.
\end{proof}

\begin{remark}
Note that in the notation of Proposition~\ref{P:H1}, 
$J$ is a one-dimensional vector space over $K_q$.
\end{remark}

\subsection{Duality and decompositions}

\begin{lemma} \label{L:irreducible}
Any irreducible $(F, \nabla)$-module over $\Gancon$
is isoclinic and quasi-constant.
\end{lemma}
\begin{proof}
An irreducible $(F, \nabla)$-module admits a slope filtration on its
underlying $F$-module by Theorem~\ref{T:slope filt}, and the slope
filtration is $\nabla$-stable by Proposition~\ref{P:isoclinic nabla}.
Hence it must have a single step, i.e., the module is isoclinic.
Since any isoclinic $(F, \nabla)$-module is quasi-constant
(Theorem~\ref{T:slope filt}), the claims follow.
\end{proof}

\begin{defn}
Let $M,N$ be $(F, \nabla)$-modules over a nearly admissible ring $S$.
By the \emph{cup product}, we will mean the natural bilinear
map $H^0_{F, \nabla}(M) \times H^1_{F,\nabla}(N) \to
H^1_{F, \nabla}(M \otimes N)$ sending
$(x, (\bv, \omega))$ to $(x \otimes \bv, x \otimes \omega)$.
Define the \emph{Poincar\'e pairing} on $M$ as the $F$-equivariant
bilinear pairing obtained by composing the cup product map
\[
H^0_{F, \nabla}(M) \times H^1_{F, \nabla}(M^\dual(-w(q)))
\to H^1_{F,\nabla}(M \otimes M^\dual(-w(q)))
\]
with the map
\[
H^1_{F,\nabla}(M \otimes M^\dual(-w(q))) \to H^1_{F,\nabla}(S(-w(q)))
\]
given by the trace map $M \otimes M^\dual \to S$.
\end{defn}

\begin{prop} \label{P:Poincare}
Assume that $k$ is algebraically closed,
and let $M$ be an $(F, \nabla)$-module over $\Gancon$.
Then the Poincar\'e pairing
$H^0_{F,\nabla}(M) \times H^1_{F, \nabla}(M^\dual(-w(q)))
\to H^1_{F, \nabla}(\Gancon(-w(q)))$ is perfect, i.e., it induces an isomorphism
\[
H^1_{F, \nabla}(M^\dual(-w(q))) \cong \Hom_K(H^0_{F,\nabla}(M),
H^1_{F,\nabla}(\Gancon(-w(q))).
\]
\end{prop}
\begin{proof}
The argument consists of a series of reductions ending with
an appeal to the calculation in
Proposition~\ref{P:H1}.
To begin with, by the snake and five lemmas, 
we may reduce to the case where $M$ is irreducible;
then $M$ is isoclinic and quasi-constant by Lemma~\ref{L:irreducible}
(note that this relies on the full theory of slope filtrations).
Let $s = c/d$ be the slope of $M'$ written in lowest terms.
Since $k$ is algebraically closed, by Theorem~\ref{T:tsuzuki},
there exists a finite separable extension
$E'$ of $E$ such that for $\Gancon' = \Gancon^{E'}$,
$([d]_* M)(-c) \otimes \Gancon'$ admits a basis of 
horizontal vectors.

Put $N = M \otimes \Gancon'$ viewed as an $(F, \nabla)$-module
over $\Gancon$; then the trace from $\Gancon'$ to $\Gancon$
induces a projector on $N$ with 
image $M$, and this map commutes with the Poincar\'e pairing.
We may thus reduce to checking the perfectness of the Poincar\'e
pairing for $N$ instead of $M$.

In other words, we have reduced to the case where $([d]_* M)(-c)$
admits a basis of horizontal vectors. At this point, we may apply
Proposition~\ref{P:change Frob} to reduce to considering a standard
Frobenius. Also, we may replace $K$ by a Galois extension (since we can
take traces down that extension); in particular, we can force $K$
to contain the $p^d$-th roots of unity.

Since $K$ contains the $p^d$-th roots of unity, we can form a trace
for the morphism $\sigma^d: \Gancon \to \Gancon$,
by averaging over automorphisms of $\Gancon$ of the form
$\{z_i\} \mapsto \zeta_i \{z_i\}$ for $\zeta_1, \dots, \zeta_m \in
\mu_{p^d}$; this gives a trace
map from $[d]^* [d]_* M$ to $M$. This means that to check perfectness
of the pairing for $M$, it is enough to do so for $[d]^* [d]_* M$.
Since the formation of $H^0$ and $H^1$ is insensitive to application
of $[d]^*$, we are reduced to checking perfectness for $[d]_* M$.

However, Theorem~\ref{T:tsuzuki} actually asserts that
$([d]_* M)(-c)$ admits a basis of horizontal vectors
\emph{stable under $F^d$}. That is, as a $(F^d, \nabla)$-module,
$[d]_* M$ splits up as a direct sum of copies of $\Gancon(c)$.
Once more by the snake lemma, we now reduce perfectness of
the Poincar\'e pairing for $[d]_* M$ to perfectness for
$\Gancon(c)$. Since the latter follows from
Proposition~\ref{P:H1}, we are done.
\end{proof}
\begin{prop} \label{P:hom ext}
Assume that $k$ is algebraically closed.
Let $M_1, M_2$ be 
irreducible $(F, \nabla)$-modules over $\Gancon$,
neither of which is isomorphic to a twist of the other.
Then 
\[
\Hom_{F,\nabla}(M_1,M_2) = \Ext^1_{F,\nabla}(M_1,M_2) = 0.
\]
\end{prop}
\begin{proof}
If $M_1$ and $M_2$ are irreducible, then $M_1$ and $M_2$ are isomorphic
if and only if $H^0_{F,\nabla}(M_1^\dual \otimes M_2) \neq 0$ if and only
if $H^0_{F,\nabla}(M_2^\dual \otimes M_1) \neq 0$.
Now Proposition~\ref{P:Poincare} yields the desired result.
\end{proof}

We may now refine the conclusion of Theorem~\ref{T:local mono} as follows.
\begin{defn}
Let $N$ be an irreducible $(F, \nabla)$-module. We say that another
$(F, \nabla)$-module $M$ is \emph{$N$-typical} if $M$ admits an
exhaustive filtration by saturated $(F, \nabla)$-submodules, 
whose successive quotients are isomorphic to twists of $N$.
If $N$ is not to be specified, we say $M$ is \emph{isotypical}.
\end{defn}

\begin{prop} \label{P:refine local}
Assume that $k$ is algebraically closed.
Let $M$ be an $(F, \nabla)$-module over $\Gancon$. Then
$M$ admits a unique direct sum decomposition
$M_1 \oplus \cdots \oplus M_l$ into isotypical
$(F, \nabla)$-submodules, such that each $M_i$ is
$N_i$-typical for some $N_i$, and no two $N_i$ are twists
of each other.
\end{prop}
\begin{proof}
The uniqueness follows from the fact that there are no nonzero
morphisms between  $(F, \nabla)$-modules which are isotypical for
irreducible modules which are not twists of each other; this follows
by repeated application of Proposition~\ref{P:hom ext}.

We prove existence by induction on $M$. 
If $M$ is irreducible, then $M$ itself is isotypical.
Otherwise, choose a short exact sequence
$0 \to M_0 \to M \to N \to 0$ with $M_0$ irreducible. Decompose
$N = \oplus N_i$ by the induction hypothesis, and let $P_i$
be the preimage of $N_i$ in $M$.
For each $i$, if $N_i$ is not $M_0$-typical,
then the exact sequence $0 \to M_0 \to P_i \to N_i \to 0$
splits, again by repeated application of Proposition~\ref{P:hom ext}.
This is true for all but possibly one $i$; we may thus split
$M$ as a direct sum of those $N_i$ plus an $M_0$-typical factor.
This completes the induction.
\end{proof}
\begin{remark} \label{R:matsuda}
One can doubtless refine Proposition~\ref{P:refine local} with more
work. For instance, it should be possible to drop the restriction that
$k$ be algebraically closed. For another, it should be possible to
show that an $M$-typical $(F, \nabla)$-module is isomorphic to
the tensor product of $M$ with a unipotent $(F, \nabla)$-module;
for true annuli, this amounts to a result of Matsuda
\cite[Theorem~7.8]{matsuda-katz}, which in turn mimics a result
of Levelt \cite{levelt} in the context of classical differential
equations.
(However, one must keep Remark~\ref{R:cannot integrate} in mind:
while Matsuda's result is purely about connections, one is compelled to
use the Frobenius also when working on fake annuli.)
This should in turn make it possible to construct a monodromy representation
in this setting, as in \cite[Theorem~4.45]{me-mono-over}, and 
perhaps to relate it to some form of the Christol-Mebkhout
construction, as in \cite[Theorem~5.23]{me-mono-over}. The latter
may be related to some conjectures of Matsuda; see
\cite{matsuda-dwork}.
\end{remark}

\subsection{Splitting exact sequences}

\begin{lemma} \label{L:rank 1}
Let
\begin{equation} \label{eq:seq1}
0 \to M_1 \to M \to M_2 \to 0
\end{equation}
be a short exact sequence 
of $F$-modules over $\Gcon[\fp]$ or $\Gancon$, and put
$d = \rank(M_1)$.
Then the sequence splits if and only if the
sequence
\begin{equation} \label{eq:seq2}
0 \to \wedge^d M_1 \to \wedge^d M \to (\wedge^d M)/(\wedge^d M_1)
\to 0
\end{equation}
splits.
\end{lemma}
\begin{proof}
If \eqref{eq:seq1} splits, then \eqref{eq:seq2} splits by the
K\"unneth decomposition. Conversely, if \eqref{eq:seq2} splits,
let $N$ be the image of $(\wedge^{d-1} M_1) \otimes M$ in $\wedge^d M$
under $\wedge$; then the exact sequence
\[
0 \to \wedge^d M_1 \to N \to (\wedge^{d-1} M_1 \otimes M_2) \to 0
\]
splits. Tensor with $M_1$ to obtain another split exact sequence
\[
0 \to (M_1 \otimes \wedge^d M_1) \to (M_1 \otimes N)
\to (M_1 \otimes \wedge^{d-1} M_1 \otimes M_2) \to 0.
\]
Twisting by $(\wedge^d M_1)^{\dual}$, we obtain a split exact sequence
\[
0 \to M_1 \to P \to (M_1 \otimes M_1^\dual \otimes M_2) \to 0
\]
for some $P$. Take the trace component within $M_1 \otimes M_1^\dual$,
tensor with $M_2$, and let $Q$ be the inverse image in $P$; we
then obtain yet another split exact sequence
\[
0 \to M_1 \to Q \to M_2 \to 0.
\]
By backtracking through the definitions, we see that this is none other
than \eqref{eq:seq1}.
\end{proof}

\begin{prop} \label{P:splitting}
Suppose that $\sigma$ is standard and that $k$ is algebraically closed.
Let $M$ be an $(F, \nabla)$-module over $\Gancon$. Then
any exact sequence
$0 \to M_1 \to M \to M_2 \to 0$ in the category of
$F$-modules splits.
\end{prop}
\begin{proof}
By Lemma~\ref{L:rank 1}, we may assume that $\rank(M_1) = 1$;
by twisting, we may assume that $M_1 \cong \Gancon$.
Also, since $k$ is algebraically closed, we may assume that
$\pi^n = q$ for some integer $n$.

Choose a coordinate system for $\Gamma$.
Let $\bv$ be an $F$-stable element of $M_1$; then
if $\bw = \Delta_1^{i_1}\cdots \Delta_m^{i_m} (\bv)$,
we have $F(\bw) = q^{-i_1 - \cdots - i_m} \bw$.
We thus obtain a nonzero
map $N \to M$ for some unipotent $(F, \nabla)$-module $N$,
whose image contains $\bv$.

Since $k$ is algebraically closed, by
Proposition~\ref{P:refine local}, we can write $M$ as a direct
sum of isotypical $(F, \nabla)$-submodules $P_1 \oplus \cdots \oplus P_l$.
At most one of the $P_i$ is isotypical for the trivial $(F,\nabla)$-module
$\Gancon$; if $P_j$ is one of the others, then the map
$N \to P_j$ obtained by composing the map $N \to M_1$ and the
projection $M \to P_j$ is zero by Proposition~\ref{P:hom ext}.
We may thus reduce to the case where $M$ is unipotent.

In this case, by repeated application of Proposition~\ref{P:H1},
we deduce that $M$ is isomorphic as an $F$-module to a direct sum
of twists of the trivial $F$-module.
Then $\bv$ must be a $K_q$-linear combination of $F$-stable generators
of summands in this decomposition; from this observation, we may
construct an $F$-stable complement
of $M_1$, yielding the desired splitting.
\end{proof}

\begin{remark}
The proof of Proposition~\ref{P:splitting} depends heavily on the hypothesis
that $\sigma$ is a standard Frobenius lift. Whether the result should even hold
otherwise is not entirely clear.
\end{remark}

\subsection{Descent of morphisms}

Following the philosophy of \cite{dejong} (as imitated in
\cite{me-full}), we now parlay our splitting
results into statements that let us descend morphisms of $(F, \nabla)$-modules
from $\Gamma$ to $\Gcon[\fp]$.

\begin{defn}
If $M$ is an $F$-module over $\Gcon[\fp]$, we may associate to $M$ two sets
of HN slopes, by passing to $\Gamma[\fp]$ and invoking
Theorem~\ref{T:slope filt} for the trivial valuation on $E$, or by
passing to $\Gancon$; we refer to these as the \emph{generic HN slopes}
and \emph{special HN slopes}, respectively. Note that the Newton
polygon of the special slopes always lies on or above that of the generic slopes,
with the same endpoint
\cite[Proposition~5.5.1]{me-slope}.
\end{defn}

\begin{prop}
Suppose that $k$ is algebraically closed.
Let $0 \to M_1 \to M \to M_2 \to 0$ be a short exact sequence
of $F$-modules over $\Gcon[\fp]$, for $\sigma$ a standard Frobenius lift,
such that $M$ acquires the
structure of an $(F, \nabla)$-module over $\Gancon$.
Suppose that each generic HN slope of $M_1$ is greater than 
each generic HN slope of $M_2$. Then the exact sequence
splits.
\end{prop}
\begin{proof}
By applying Lemma~\ref{L:rank 1} and taking duals, 
we may assume that $\rank(M_2) = 1$;
by twisting, we may assume that $M_2 \cong \Gcon[\fp]$. Then
the slopes condition is that $M_1$ has all generic HN slopes positive.
By Proposition~\ref{P:splitting}, the exact sequence
\[
0 \to M_1 \otimes \Gancon \to M \otimes \Gancon \to \Gancon \to 0
\]
of $F$-modules splits; by \cite[Proposition~7.4.2]{me-slope}, 
the original sequence also splits.
\end{proof}

\begin{theorem} \label{T:full faithful}
\begin{enumerate}
\item[(a)]
Let $M$ be an $F$-module over $\Gcon[\fp]$, for $\sigma$ a standard
Frobenius lift, which acquires the structure
of an $(F, \nabla)$-module over $\Gancon$.
Then the natural map
\[
H^0_{F}(M) \to H^0_{F}(M \otimes \Gamma[\fp])
\]
is a bijection.
\item[(b)]
Let $M$ be an $(F, \nabla)$-module over $\Gcon[\pi^{-1}]$.
Then the natural map
\[
H^0_{F,\nabla}(M) \to H^0_{F,\nabla}(M \otimes \Gamma[\fp])
\]
is a bijection.
\end{enumerate}
\end{theorem}
\begin{proof}
In either case, there is no harm in assuming that $k$ is algebraically closed.
\begin{enumerate}
\item[(a)]
Given $\bv \in H^0_F(M \otimes \Gamma[\fp])$, we obtain
an $F$-equivariant dual map $\phi: M^\dual \to \Gamma[\fp]$.
Let $N_0$ be the kernel of $\phi$; by
\cite[Proposition~7.5.1]{me-slope} (applicable because any monomial
field admits a valuation $p$-basis, namely any coordinate system), the preimage
$N_1 = \phi^{-1}(\Gcon[\fp])$ has the property that
$N_1/N_0 \cong \Gcon[\fp]$, and $M^\dual/N_1$ has all
generic slopes negative.

By Proposition~\ref{P:splitting}, $N_1$ admits an $F$-stable
complement in $M^\dual$, so $N_1/N_0$ admits an $F$-stable
complement $P$ in $M^\dual/N_0$.
However, the generic slopes of $P$ are the same as those of $M^\dual/N_1$, so
they are all negative. Thus the map
$P \to \Gamma[\fp]$ obtained by composition
with $\phi$ is forced to vanish by \cite[Proposition~7.5.1]{me-slope},
whereas $\phi: M^\dual/N_0 \to \Gamma[\fp]$ is injective.
We must then have $P=0$ and $N_1 = M^\dual$, so
$\phi$ maps into $\Gcon[\fp]$ and $\bv \in M$, as desired.
\item[(b)]
By Proposition~\ref{P:change Frob}, there is no harm in reducing
to the case of $\sigma$ standard. The result now follows from (a).
\end{enumerate}
\end{proof}

\begin{remark}
Note that even for true annuli, 
Theorem~\ref{T:full faithful} makes an assertion (namely (a)) not
covered by \cite{me-full}.
\end{remark}

\end{document}